\documentclass[11pt,a4paper,leqno]{amsart}

\usepackage{fullpage,amsmath,amsthm,amsfonts,amscd,
latexsym,amssymb,eucal,curves}
  \usepackage[all]{xy}
\usepackage[dvips]{graphicx}
\usepackage[dvips]{color}
\newcommand{\ZZ}{\mathbb{Z}}

\newcommand{\CC}{\mathbb{C}}

\newcommand{\PP}{\mathbb{P}}

\newcommand{\HH}{\mathbb{H}}

\newcommand{\LL}{\mathbb{L}}

\newcommand{\MM}{\mathbb{M}}

\newcommand{\QQ}{\mathbb{Q}}
\newcommand{\RR}{\mathbb{R}}

\newcommand{\VV}{\mathbb{V}}
\newcommand{\Ker}{{\rm Ker}}

\newcommand{\End}{{\rm End}}

\newcommand{\id}{{\rm id}}
\newcommand{\dR}{{\rm dR}}

\newcommand{\Gal}{{\rm Gal}}

\newcommand{\Sym}{{\rm Sym}}

\newcommand{\fraco}{\mathfrak{o}}

\newcommand{\cLL}{{\mathcal L}}
\newcommand{\cOO}{{\mathcal O}}
\newcommand{\cEE}{{\mathcal E}}
\newcommand{\cFF}{{\mathcal F}}

\newcommand{\SL}{{\rm SL}}

\newcommand{\diag}{{\rm diag}}
\newcommand{\mfC}{{\mathfrak C}}

\newcommand{\Sp}{{\rm Sp}}

\newcommand{\rel}{{\rm rel}}
\newcommand{\tr}{{\rm tr}}
\newcommand{\sms}{\smallsetminus}

\newcommand{\codim}{{\rm codim}}

\newcommand{\ve}{{\varepsilon}}
\newcommand{\ol}{\overline}

\newcommand{\GL}{{\rm GL}}

\renewcommand{\div}{{\rm div}}
\newcommand{\red}{{\rm red}}
\newtheorem{Defi}{Definition}[section]

\newtheorem{Rem}[Defi]{Remark}
\newtheorem{Example}[Defi]{Example}
\newtheorem{Prop}[Defi]{Proposition}
\newtheorem{Lemma}[Defi]{Lemma}
\newtheorem{Cor}[Defi]{Corollary}
\newtheorem{Thm}[Defi]{Theorem}

\newtheorem{Q}[Defi]{Question}

\makeatletter
\renewcommand{\subsection}{\@startsection{subsection}{2}%
        {\z@}{-3.25ex plus -1ex minus-.2ex}{-1em}{\bf}}
\makeatother

\begin{document}{\large}
\title{Linear manifolds
in the moduli space of one-forms}
\author{Martin M\"oller}

\begin{abstract}
We study closures of $\GL^+_2(\RR)$-orbits in the
total space $\Omega M_g$ of the Hodge bundle over the
moduli space of curves under the assumption that they
are algebraic manifolds.
\par
We show that, in the generic stratum,
such manifolds are the whole stratum, the hyperelliptic
locus or parameterize curves whose Jacobian has additional
endomorphisms. This follows from a cohomological description of the 
tangent bundle to $\Omega M_g$. For non-generic strata
similar results can be shown by a case-by-case inspection.
\par
We also propose to study a notion of 'linear manifold' that comprises
Teichm\"uller curves, Hilbert modular surfaces and the
ball quotients of Deligne and Mostow. Moreover, we give an 
explanation for the difference between Hilbert modular surfaces 
and Hilbert modular threefolds with respect to this notion
of linearity.
\date{\today}
\end{abstract}
\maketitle

\section*{Introduction}
The total space $\Omega M_g$ of the Hodge bundle over $M_g$ or equivalently
the space of pairs $(X,\omega)$ of a Riemann surface $X$ and
a holomorphic one-form $\omega \in \Gamma(X,\Omega^1_X)$
admits a linear structure. That is, integration of one-forms
along a basis of the first homology relative to zeroes of $\omega$
gives locally a map to from $\Omega M_g$ to $\CC^N$ and different
choices of the basis yields a linear coordinate change, in fact defined over 
$\RR$. We review in Section~\ref{affstruc} equivalent ways of interpreting 
a linear structure in terms of connections or local systems.
\par
Due to the linear coordinate changes defined over $\RR$, 
there is a natural action of $\GL_2^+(\RR)$
on $\Omega M_g$, or rather the complement $\Omega M_g^*$ of the zero
section. A guiding question consists in the classification of
the closures of these orbits.
The situation seems to have several similarities to Ratner's theorem
on orbit closures of homogeneous manifolds. Indeed for $g=2$
the classification has been achieved by McMullen (\cite{McM03b}).
\par
Here, to analyze what happens for $g \geq 3$,  we propose to split the 
orbit closure problem into three subquestions. First, to show that 
orbit closures are indeed complex manifolds. Second, to show that 
these complex manifolds are in fact algebraic. Third, to classify
these algebraic manifolds. As Kontsevich observed,
orbit closures that are complex manifolds are 
submanifolds of $\Omega M_g$, or more precisely of its strata,
that inherit a linear structure defined over $\RR$.
See Section~\ref{affstruc} for details. The converse also holds:
If a submanifold of a stratum inherits a linear structure defined over $\RR$,
it is $\GL_2^+(\RR)$-invariant.  Thus, the third problem
may be translated into a purely algebro-geometric problem, 
interesting independently of the orbit closure question. 
The manifolds with linear structure include both Hilbert modular
surfaces (\cite{McM03b}) and, as we will show, ball quotients (\cite{DeMo86}).
\par
The purpose of the present paper is to stress the role of
the hyperelliptic locus in this classification problem. One
of our main results is:
\par
{\bf Theorem~\ref{comphypend}.}
Suppose that $B$ is an algebraic submanifold of the generic
stratum $S=\Omega M_g(1,\ldots,1)$ with linear structure. Then
there are three possibilities:
\begin{itemize}
\item[(i)] $B$ is a connected component of $S$.
\item[(ii)]  $g \geq 3$ and $B$ is the preimage in $S$
of the hyperelliptic locus in $M_g$. 
\item[(iii)] $B$ parameterizes curves with a Jacobian whose endomorphism ring
is strictly larger than $\ZZ$.
\end{itemize}
\par
For non-generic strata a similar statement holds true, maybe up
to some exceptional linear submanifolds in strata whose projection
to $M_g$ has small fiber dimension. We do a case by case analysis
for $g=3$ in Section~\ref{secHypEnd}.
\par
In order to prove Theorem~\ref{comphypend} we describe in 
Section~\ref{tangentbundle}
the tangent bundle cohomologically. This has been done pointwise
already in \cite{HuMa79} and yields, locally, the period coordinates
used throughout when dealing with the $\GL_2^+(\RR)$-action.
The advantage of the global viewpoint consists of keeping track
of the twist by $\cOO(1)$ on $\Omega M_g$, or rather on a
projective completion. We thus obtain some information
on the foliation of $\Omega M_g$ by constant absolute
periods via the multiplication map of one-forms.
From this viewpoint, the hyperelliptic locus figures among
the above list, precisely because the multiplication map
on one-forms fails to be surjective there for $g \geq 3$.
\par
This paper grew out of the attempt to understand the difference,
discovered in \cite{McM03a},
between the linearity of (the eigenform bundle) over a Hilbert 
modular surface and the non-linearity of the Hilbert modular threefold
parameterizing abelian varieties of dimension three with real
multiplication by $\ZZ[\zeta_7+\zeta_7^{-1}]$. 
\par
For Hilbert modular threefolds we show that there are
two obstacles to linearity (see Theorem~\ref{obstacles}). The
first one is related to properties of the multiplication map.
It is already used in \cite{McM03a}. We show that that this
phenomenon can arise only at the hyperelliptic locus. Geometrically
it says that the normal vector to the hyperelliptic locus coincides
with one of the three natural foliations.
\par
The second one is related to the intersection of the eigenlocus
over Hilbert modular threefold with the foliation by constant absolute
periods. Neither of the problems appears in genus $g=2$ and this
observation enables us to reprove linearity in genus $g=2$
(Proposition~\ref{HM2linear}, \cite{McM03a} Theorem~7.1).
\par
It seems very likely that the generic part of a Hilbert modular
threefold is never linear. But a proof of this might
first need to be able to decide in which stratum a Hilbert modular
threefold lies generically, depending on the endomorphism ring.
At present, this question remains open.
\par
In Section~\ref{notoverRR} we examine which manifolds we get
if we replace linearity defined over $\RR$ by just linearity. 
As remarked above,
the ball quotients of Deligne and Mostow fit into this picture.
But contrary to the case of linear transition functions defined 
over $\RR$, pathological cases and
manifolds with linear structure unrelated to uniformization
occur as well. We propose a definition of a linear manifold that
includes precisely those manifolds with a linear structure
related to the uniformization of the manifold.
The condition we demand additionally is the existence of a
compactification to which the linear structure ``extends''
in terms of a surjection of some tangent bundle to a Deligne
extension of a local system.  
\par
Already in the case of $2$-dimensional linear manifolds the
definition has surprising consequences. Using the cohomological 
description we quickly reprove in Section~\ref{lowg2} the classification
result of McMullen in genus $2$, of course under the general
assumption of algebraicity. But there are, even in genus $2$, two-dimensional
linear manifolds that are not canonical lifts of Teichm\"uller curves, since
the linear structure is not defined over $\RR$.
\par
In Section~\ref{lininH22} we apply the same classification techniques to
the hyperelliptic locus $\mathcal{H}$ inside the odd spin component of 
the stratum $\Omega M_3(2,2)$. This locus is maybe the
most simple, besides $g=2$, to study $\GL_2^+(\RR)$-orbit closures.
See also \cite{HLM06}. One consequence is:
\par
{\bf Corollary \ref{pAclosure}.}
Let $\Delta$ be a Teichm\"uller disc 
generated by a pair $(X_0,\omega_0)  \in \mathcal{H}$ which
is stabilized by a pseudo-Anosov diffeomorphism with trace field
of degree $3$ over $\QQ$. If the closure of the orbit 
$\GL_2^+(\RR) \cdot (X_0,\omega_0)$ is an algebraic manifold,
then it is either the canonical lift of a Teichm\"uller curve to $\mathcal{H}$
or it is as big as possible, i.e.\ 
$$\ol{\GL_2^+(\RR) \cdot (X_0,\omega_0)} = 
\mathcal{H}.$$
\par
{\bf Acknowledgements.} The author had many fruitful discussions with 
Eckart Viehweg on linear manifolds and Hilbert modular varieties. 
The author thanks him and H\'el\`ene Esnault a lot for their help and the 
participants of the conference 'Outer space and Teichm\"uller space', 
CIRM (2007), and Curt McMullen for their comments. The author also thanks
the referees for their comments.
\par
\section{Manifolds with affine structures and linear structures} 
\label{affstruc}

{\bf Affine and linear structures}
\newline
An {\em affine structure} on a complex manifold $M$ (say 
of dimension $d$) is an atlas of charts $(U_i, \varphi_i: U_i \to \CC^{d})$
all whose transition functions $\varphi_j\circ \varphi_i^{-1}$
are affine maps. For convenience we suppose that
the atlas is maximal for that property.
\par
We recall below that it is equivalent to have a flat connection 
$\nabla: \Omega_M \to \Omega_M \otimes \Omega_M $ on $M$
which is symmetric, i.e.\ sends closed forms to ${\rm Sym}^2 \Omega_M$.
Another equivalent definition is that there is a $\CC$-local 
system $\LL$ on $M$ such that $T_M = \LL \otimes_\CC \cOO_M$ 
which is torsion-free, i.e.\ such that the corresponding 
connection $\nabla^{T_M}$ satisfies
$$\nabla^{T_M}_X Y - \nabla^{T_M}_Y X = [X,Y]$$
for any local sections $X,Y$ of $T_M$.
\par
To see that the definitions are equivalent, start with a manifold
$M$ with an affine structure. The $\CC$-valued functions that
are locally affine form a local system ${\rm Aff}_M$ of rank $d+1$ 
that contains the constant functions. The quotient ${\rm Aff}_M/\CC$ is
a local system whose associated vector bundle is the cotangent bundle.
The symmetry of this connection is checked locally.
\par
Conversely starting with a connection $\nabla$, consider
the functions whose total differential is flat for $\nabla$. 
Locally, a complement of the constant functions defines a
chart. These charts fit to an affine structure if the
connection is symmetric.
\par
The equivalence between $\nabla$ being symmetric and $\nabla^{T_M}$
being torsion free is immediate from the definition of the dual connection.
We will drop the superscript $T_M$ from the connection $\nabla$ in the sequel.
\par
A {\em linear structure} (sometimes also called {\em radiant
affine structure}) is defined by an atlas as above, whose
transition functions are linear. This is equivalent to the
monodromy group of $\LL$ being equal to the {\em holonomy group}, the
monodromy group of ${\rm Aff}_M$. Again equivalently, $\CC$
has a complement in ${\rm Aff}_M$.
A manifold $M$ has a {\em linear structure defined over
$\RR$} or {\em
is locally defined by $\RR$-linear equations}
if the transition maps are in $\GL_d(\RR) \subset \GL_d(\CC)$.
\par
We will reserve the expression 'linear (sub)manifold' for a manifold
with linear structure plus additional conditions related to the boundary,
see Section~\ref{notoverRR}. 
\par
\begin{Lemma} \label{critlinsubmf}
If $B \subset M$ is a submanifold of an affine manifold $M$
and the connection $\nabla$ restricts to a connection
$\nabla_B: \Omega_B \to \Omega_B \otimes \Omega_B$, 
then $B$ has also an affine structure. If moreover $M$ has
a linear structure, then $B$ has a linear structure.
\end{Lemma}
\par
\begin{proof} Since $\nabla_B$ is obtained by restriction from
$\nabla$, it is automatically flat and torsion free. The
second claim is obvious.
\end{proof}
\par
{\bf Infinitesimal $\GL_2^+(\RR)$-action}
\newline
A complex manifold $M$ with a linear 
structure defined over $\RR$ automatically has a natural infinitesimal
$\GL_2^+(\RR)$-action. I.e.\ at each point $x \in M$ there
is a neighborhood of ${\rm Id} \in \GL_2^+(\RR)$ that acts
on a neighborhood of $x$ in the following way:  On
the linear charts $\varphi_i(x) \in \CC^d = \RR^d \otimes_\RR \CC$  
the group $\GL_2^+(\RR)$ acts linearly on $\CC \cong \RR^2$. This
action commutes with the transition functions, since they
are defined over $\RR$, they act on the first factor of the tensor
product and since $\RR$, diagonally embedded in $\GL_d(\RR)$
resp.\ $\GL_2^+(\RR)$, is central in both cases. This argument
also implies that the action, if defined, does not depend on the chart
chosen.
\newline
One cannot hope to globalize this action without further knowledge
on $M$. For example one could have removed a closed subset, but
not an orbit, from a manifold where the global action was defined.
\par
{\bf The Hodge bundle over $M_g$ and its strata.}
\newline
Let $\Omega M_g$ denote (the total space of) the Hodge bundle over the
moduli space $M_g$ of curves of genus $g$ and let $\Omega M_g^*$ 
denote the complement of the zero section. Points in 
$\Omega M_g^*$ correspond to pairs $(X,\omega)$ of 
a curve $X$ of genus $g$ and a non-zero holomorphic one-form on $X$.
\newline
There is a natural $\GL_2^+(\RR)$-action on $\Omega M_g^*$
defined by post-composing the local charts given by
integrating $\omega$ with the linear map. See \cite{MaTa02} or
\cite{Zo06} for recent surveys.
This action respects the strata  $\Omega M_g(k_1,\ldots,k_r)$
consisting of pairs $(X,\omega)$ such that the zeroes of $\omega$
are of the form $\sum_{i=1}^r k_i P_i$ for disjoint points $P_i$
on $X$. The connected components of the strata have
been completely determined in \cite{KoZo03}.
\par
{\bf Period map.} 
\newline
Let $S_{T_g} \subset \Omega T_g$ be the pullback of a stratum $S$
the universal bundle to the Teichm\"uller space $T_g$. On a
simply connected open subset $U$ of $S_{T_g}$ choose a basis of the relative 
homology. Integrating the one-from along this
basis gives a map, called the {\em period map}, that attaches to a
pair $(X,\omega)$ a point in $\CC^N$ for $N=2g-1+r$. 
It is well-known that the period map is a local biholomorphism
(\cite{Ve86}). One may deduce this also from the considerations
in Section~\ref{tangentbundle}. The mapping class group
acts linearly on the relative periods. Consequently,  the period map provides
the stratum with a linear structure. 
\par
The mapping class group maps absolute periods to absolute periods.
Hence, forgetting relative periods gives locally a map 
$\CC^N \to \CC^{2g}$, equivariant with respect to 
the action of the mapping class group. Thus, the strata
have a natural foliation by leaves of complex dimension $r-1$,
called the {\em foliation by constant absolute periods}.
\par
Our starting point is the following observation. Folklore attributes
the application of it to the period coordinates of
a $\GL_2^+(\RR)$-invariant submanifold to Kontsevich:
\par
\begin{Prop} \label{Kontsevich}
A $\GL_2^+(\RR)$-invariant $B$ closed analytic subspace
in $\CC^N$ is linear and defined locally by $\RR$-linear equations.
\end{Prop}
\par
{\bf Proof:}
Suppose that $B \subset \CC^N$ is a hypersurface given 
by a power series $f$ with complex coefficients in the variables
$z_1=x_1+iy_1,\ldots,z_N=x_N+iy_N$. We claim that $f$ is a 
homogeneous  polynomial of some degree $m$. To show this assume the
contrary. Then for some $\lambda \neq 1$, the power series
$f$ and $f(\lambda z_1,\ldots,\lambda z_N)$ are linearly independent.
This contradicts that $B$ is a hypersurface and invariant under
the diagonal subgroup in $\GL_2^+(\RR)$.
\par
The natural action of $\GL_2^+(\RR)$ on points of $\CC^N\cong\RR^N
\otimes_\RR \RR^2$ by the linear action on the second factor corresponds on 
the level of polynomial defining closed subspaces to the 
action 
$$ g \cdot f := f \circ g^{-1} := f((\tilde{a}x_1 + \tilde{b}y_1) +
i (\tilde{c}x_1 + \tilde{d}y_1),\ldots, (\tilde{a}x_N + \tilde{b}y_N) +
i (\tilde{c}x_N + \tilde{d}y_N)), $$
where $g^{-1} = \left(\begin{matrix} \tilde{a} & \tilde{b} \\
\tilde{c} & \tilde{d} \\ \end{matrix} \right)$.
Since $f$ is a homogeneous polynomial of degree $m$, the $\RR$-vector space 
$$V = \left\langle g\cdot f, \,g \in \GL_2^+(\RR)\right\rangle$$ 
of polynomials 
with complex coefficients contains a vector of weight $m$ 
(as a $\GL_2^+(\RR)$-module).
Hence $\dim_\RR V \geq m+1$. We claim that this implies
that that the real codimension
of $B$ in $\CC^N$ is at least $m+1$. Since $B$ is
a complex hypersurface we conclude that $f$ is linear
and moreover that we can choose $f$ to have real coefficients.
\par
To prove the claim, suppose that the zero set in $\RR^{2N}$ of 
the polynomials in $V$ in the variables $x_i$ and $y_i$ has codimension 
less than $m+1$. Tensoring with $\CC$, we conclude that the zero 
set in $\CC^{2N}$ of the polynomials in $V$ has codimension 
less than $m+1$ and we obtain a contradiction.
\par
We use the hypersurface case to argue by induction on both
$N$ and $\codim_{\CC^N}(B)$. The map $\varphi: \CC^N \to
\CC^{N_1}$ by forgetting the last coordinate is
$\GL_2^+(\RR)$-equivariant. Hence the image $\varphi(B)$ has a
linear structure defined over $\RR$ by
induction hypothesis. If we let $W$ be the preimage $\varphi(B)$
then $W\cong \CC^{N'}$ for some $N'$. Moreover, the 
embedding $\CC^{N'} \to \CC^N$ is given by $\RR$-linear equations.
On the other hand, $B$ either equals $W$ or a hypersurface
in $W$ and we can again apply the induction  hypothesis. 
\hfill $\Box$
\par
We emphasize that for a submanifold $B$ in a manifold $S$ with linear
structure (usually a stratum) we say {\em $B$ has a linear
structure} as shorthand for '$B$ inherits a linear structure from $S$'
or '$B$ is defined locally by linear equations in period coordinates'.
That is, we are never interested in intrinsic linear structures on $B$.
\par
The converse of Proposition~\ref{Kontsevich} holds in the 
strata of $\Omega M_g$.
\par
\begin{Prop} \label{Kontconverse}
Let $B \subset \Omega M_g(k_1,\ldots,k_r)$ be a closed, linear
submanifold in some stratum. If $B$ is defined by $\RR$-linear
equations then $B$ is $\GL_2^+(\RR)$-invariant.
\end{Prop}
\par
\begin{proof} On $\Omega M_g(k_1,\ldots,k_r)$ we have an action
of  $\GL_2^+(\RR)$, not only of a neighborhood of ${\rm Id} \in 
\GL_2^+(\RR)$. By the
hypothesis this action respects $B$. 
\end{proof}
\par
\begin{Rem} {\rm
We a priori work in the category of complex analytic manifolds.
Some of the stronger classification results below only work
under the hypothesis that these manifolds are algebraic.
Dividing by the $\CC^*$-action on $\Omega M_g^*$ one observes
that these two kinds of objects coincide, if one can show
that the closures of objects in question (linear manifolds in strata)
in some projective compactification are still manifolds.
}\end{Rem}
\par
{\bf Orbifolds.}
\newline
The moduli space of curves $M_g$ should be considered as a complex 
orbifold in the sequel. To avoid technicalities and since all our results
will be independent of passing to a finite unramified cover
of $M_g$ we fix a suitable one (say $M_g^{[n]}$ with some 
level-$[n]$-structure),
which is a smooth manifold and over which the universal
family of curves $f: X \to M_g$ exists. We nevertheless
keep the symbol $M_g$ for simplicity. We also denote the pullback of
the  universal family to any manifold over $M_g$ by
$f$, hoping not to create confusion.
\par
However with this notation the Torelli map $t^{[n]}: M^{[n]}_g \to A^{[n]}_g$
is no longer injective. Instead it is a $2:1$ covering ramified
precisely over the hyperelliptic locus, see diagram \eqref{levelstructure}.
\par
{\bf Variation of Hodge structures}
\newline
We will mainly be concerned with variations of Hodge structures (VHS)
of weight one on
an algebraic manifold $B$, sometimes completed by a normal
crossing divisor $S$ to $\ol{B}$.
A weight one VHS consists of a $\ZZ$-local system $\VV$ plus a filtration
of vector bundles $\VV^{(1,0)} \subset (\VV \otimes \cOO_B)_{\rm ext}=:
{\mathcal V}$, 
where ${\rm ext}$ denotes the Deligne extension (\cite{De70}). The only compatibility
condition of the filtration we make use of here is 
$\VV^{(0,1)}:= {\mathcal V}/ \VV^{(1,0)} \cong (\VV^{(1,0)})^\vee$.
With this data one defines the Higgs field to be the composition
$$\theta: \VV^{(1,0)} \to {\mathcal V}
\to {\mathcal V} \otimes \Omega^1_B(\log S)
\to {\mathcal V}/ \VV^{(1,0)}  \otimes \Omega^1_B(\log S).$$ 
If the VHS comes from a family of smooth curves $f:X \to B$, i.e.\
if $\VV = R^1 f_* \ZZ$ and $\VV^{(1,0)} = f_* \omega_{X/B}$, we can
use duality to obtain from $\theta$ a map 
$$ {\rm Sym}^{2}(\VV^{(1,0)}) \to \Omega^1_B(\log S),$$
which is well-known to factor through $f_* \omega_{X/B}^2$. The dual
of this factorization is usually called the {\em Kodaira-Spencer map.}
\par
\section{Cohomological description of several 
tangent bundles} \label{tangentbundle}

The goal of this technical section is to describe the (co)tangent bundle 
to strata in $\Omega M_g$ by hypercohomology sheaves. It will be
important to do this in the relative setting not just pointwise.
Only in that way we can keep track of the twist by $\cOO_\PP(1)$.
Although this bundle admits a trivialization over the strata, 
this twist helps to identify certain bundle maps as induced by
multiplication maps on one-forms. Properties of the multiplication 
map will then be exploited in the subsequent sections.
\par

\subsection{Tangent bundle to the completed one-form bundle}

Our first aim is to describe the tangent bundle to $\Omega M_g$.
For technical reasons we complete this bundle to
a projective bundle, rather than taking the projectivisation 
$\PP \Omega M_g$. The reason is that the latter has no longer
a linear structure. We start fixing some conventions.
\par
For a vector bundle $\cEE$ on $Y$ let $\PP$ be the
projective bundle ${\rm Proj}({\rm Sym(\cEE^\vee)})$
with projection $p: \PP \to Y$.
There is canonical map $ \cOO_{\PP}(-1) \to p^* \cEE$
and the Euler exact sequence 
$$ 0 \to \cOO_{\PP} \to p^* \cEE(1) \to T_{\PP/Y} \to 0.$$
Suppose that $\cEE = \cEE' \oplus \cOO_Y$. Then  the 
hyperplane at infinity is defined by applying ${\rm Proj}$
to the inclusion $ \cOO_Y \to  \cEE' \oplus \cOO_Y$ on the
second summand. The exact sequence
\begin{equation} \label{logrem}
 0 \to \Omega^1_{\PP/Y} \to \Omega^1_{\PP/Y}(\log \infty) \to \cOO_\infty
\to 0
\end{equation}
together with the Euler sequence implies that
$\Omega^1_{\PP/Y}(\log \infty) \cong p^* \cEE'^\vee (-1).$
\par
We apply the above remarks in the case $Y= M_g$ and $\cEE = 
(f_* \Omega^1_{X/M_g} \oplus \cOO_{M_g})$, where $f: X \to M_g$
is the universal family of curves. Moreover we abuse the letter $f$
also for pullback families, e.g. $f: X_\PP \to \PP$.
\par
\begin{Thm} \label{tangtoPP}
The tangent bundle of $\PP = 
\PP((f_* \Omega^1_{X/M_g} \oplus \cOO_{M_g})^\vee)$ sits in an exact sequence 
$$ 0 \to T_{\PP/{M_g}}(-\log \infty) \to 
T_{\PP}(-\log \infty) \to p^* T_{M_g} \to 0.$$ 
The extension comes form 
$$ T_{\PP} = R^1 f_* (\varphi: (\Omega^1_{X_\PP/\PP})^\vee 
\to f^* \cOO_{\PP}(1) \otimes  \Omega_{X_\PP/\PP}^1) $$ 
and the identifications
$$T_{\PP/{M_g}}(-\log \infty) = f_* \Omega^1_{X_\PP/\PP}(1)
\; \text{and}\;  p^* T_{M_g} = R^1 f_* (\Omega^1_{X_\PP/\PP})^\vee $$
The map $\varphi$ is given by 
$$\xymatrix{
(\Omega^1_{X_\PP/\PP})^\vee  \ar[r]^{{\rm can}} & f^* \cOO_{\PP}(1) \ar[rr]^(.3){d}
&& f^* \cOO_{\PP}(1) \otimes \Omega_{X_\PP}^1(\log f^{-1}(\infty)) \ar[r]
& f^* \cOO_{\PP}(1) \otimes \Omega_{X_\PP/\PP}^1}.
$$ 
\end{Thm}
\par
\begin{proof} We first show that the claims are true locally
at some point $(X^0,\omega^0) \in \PP(\CC)$.
In this part we basically follow \cite{HuMa79}, who deal with the case of
quadratic differentials. See also \cite{We83}.
\newline
Let $X_\ve$ be a deformation of $X^0$ over $\Delta := {\rm Spec}\;
\CC[\ve]/\ve^2$ and let $\omega_\ve \in \Gamma(X_\ve,\Omega^1_{X_\ve/\Delta})$ 
be a deformation of $\omega^0$. We want to associate with these data
a class in $H^1(X_0,T_{X^0} \to \Omega^1_{X_0})$, where the map
is the Lie derivative, i.e.\ contracts a tangent vector against
$\omega^0$ and then takes the exterior derivative of the resulting function.
\newline
Let $\{U_\alpha ={\rm Spec}\;  A_\alpha \times \Delta \}$ be a covering
of $X_\ve$ such that $X_\ve$ is given by transition functions
$\psi_{\alpha \beta}$ with
$$ \psi_{\alpha \beta}(\ve) = \ve, \quad
\psi_{\alpha \beta}(f) = f + \ve D_{\alpha \beta}(f). $$
Since $\psi_{\alpha \beta}$ are ring homomorphisms, 
$D_{\alpha \beta}$ is a $\CC$-derivation with values in $A_{\alpha \beta}$.
The transition functions of $\Omega^1_{X_\ve/\Delta}$ are given by
$$ df + \ve dg \mapsto d(f + \ve D_{\alpha \beta}(f)) + \ve d
(g + \ve D_{\alpha \beta}(g)) = df + \ve dg + \ve d D_{\alpha \beta}(f).  $$  
If we describe the section as 
$(\omega_\ve)_\alpha = \omega_\alpha^0 + \ve ds_\alpha$, the gluing
condition of the bundle is
$$ ds_\alpha - ds_\beta = d D_{\alpha \beta}(\omega^0_\alpha).$$
Hence the pair $(\{D_{\alpha \beta} \}, \{ds_\alpha\})
\in C^1(U_\alpha,T_X) \oplus C^0(U_\alpha, \Omega^1_{X^0})$ forms
a $1$-cochain of the desired complex. It is straightforward to
see that the pair is indeed a cochain and that different choices
modify the cochain only up to a coboundary (see loc.\ cit.). 
One obtains in this way an identification of the tangent space
$$T_{(X^0,\omega^0)} \PP \cong H^1(X_0,T_{X^0} \to \Omega^1_{X_0}).$$
\par
What remains to prove the proposition is to show that this identification
fits together in the bundle. In order to do so we define $\varphi$ as 
the composition map
$$\xymatrix{
(\Omega^1_{X_\PP/\PP})^\vee  \ar[r]^{{\rm can}} & f^* \cOO_{\PP}(1) \ar[rr]^(.3){d}
&& f^* \cOO_{\PP}(1) \otimes \Omega_{X_\PP}^1(\log f^{-1}(\infty)) \ar[r]
& f^* \cOO_{\PP}(1) \otimes \Omega_{X_\PP/\PP}^1}
$$ 
and consider it as a complex $\mfC^\bullet$ in degrees $0$ and $1$.
Here ``${\rm can}$'' is the contraction against the universal section, 
i.e, the dual of 
$$ f^* \cOO_{\PP}(-1) \to f^*( p^* (f_* \Omega^1_{X/M_g} \oplus \cOO_{M_g}))
\to f^* (p^* f_* \Omega^1_{X/M_g}) =  f^* f_* \Omega^1_{X_\PP/\PP} 
\to \Omega^1_{X_\PP/\PP} $$ 
and $d$ is the exterior derivative.
\par
We will define a map $T_\PP \to R^1 f_* (\mfC^\bullet)$
that boils down to the one described by deformation
theory for any thick point $\Delta \to \PP$. All we need is to construct
a map $\widetilde{\varphi}$, that fits into the following exact
sequence of complexes
$$ \xymatrix{ 0 \ar[r] & T_{X_\PP/\PP} \ar[r] \ar[d]^\varphi & T_{{X_\PP},f}
(-\log f^{-1}(\infty)) \ar[rr] \ar[d]^{\widetilde{\varphi}} 
&& f^{-1} T_\PP(-\log \infty) \ar[r] \ar[d] & 0 \\
0 \ar[r] &  f^* \cOO_{\PP}(1) \otimes \Omega_{X_\PP/\PP}^1 \ar@2{-}[r] &  
f^* \cOO_{\PP}(1) \otimes \Omega_{X_\PP/\PP}^1 \ar[rr] && 0 \ar[r] & 0 \\
}$$
where the upper line is the restriction of the standard exact sequence
from $f^* T_\PP(-\log \infty)$ to $f^{-1} T_\PP(-\log \infty)$.
Once we have established this diagram the coboundary map associated
to $R f_*(\cdot)$ will provide the desired map. Observe that
after applying $R f_*(\cdot)$ there is no difference between
$f^* T_\PP(-\log \infty)$ and $f^{-1} T_\PP(-\log \infty)$, but
to construct ${\widetilde{\varphi}}$ the difference is essential.
\par
For any small enough open set $U$  elements of 
$T_{{X_\PP},f} (-\log f^{-1}(\infty))$ are of the form
$\lambda t + \tau$ for $\tau \in f^{-1}T_\PP(-\log \infty)(U)$, $t \in
T_{X_\PP/\PP}(U)$ and $\lambda \in \cOO_{X_\PP}(U)$. Contraction
of $\lambda t + \tau$ against an element of $\Omega^1_{X_\PP/\PP}$
gives a function on $U$, well-defined only up to $f^{-1}$ of functions
on $\PP$. Hence $d \circ {\rm can}$ is well-defined up
to a section of 
$f^* \cOO_{\PP}(1) \otimes f^{-1}(\Omega^1_\PP(\log \infty))(U)$.
Consequently, the composition with the projection 
to $f^* \cOO_{\PP}(1) \otimes \Omega_{X_\PP/\PP}^1$ is indeed the
well-defined map $\tilde{\varphi}$ we need.
\par
A two-term complex $\cEE_0 \to \cEE_1$ in degrees zero and one admits the
so-called stupid filtration by subcomplexes
$$ \{0\} \subset [0 \to \cEE_1] \subset [\cEE_0 \to \cEE_1].$$
This filtration naturally defines a short exact sequence of complexes.
Applied to $\mfC^\bullet$ we obtain the commutative diagram
$$ \xymatrix{
 0 \ar[r]&  T_{\PP/M_g}(-\log \infty) \ar[r] \ar[d] & T_\PP
\ar[r] \ar[d] &  p^* T_{M_g} \ar[r]  \ar[d] & 0 \\
 0 \ar[r]&  f_* \Omega^1_{X_\PP/\PP} (1) \ar[r] & R^1 f_* (\mfC^\bullet) 
\ar[r] & R^1 f_* (\Omega^1_{X_\PP/\PP})^\vee  \ar[r] & 0. \\
} $$
By deformation theory, as above simply forgetting $\omega_\ve$,
the right vertical arrow is an isomorphism too. This proves
the remaining statements claimed in the theorem.
\end{proof}
\par

\subsection{Tangent bundle to a stratum}
\label{tangstrat}

Let now $S:=\Omega M_g(k_1,\ldots,k_r) \hookrightarrow \Omega M_g 
\hookrightarrow \PP$ be a connected component of a stratum. 
The dimension of $S$ is
$2g+r-1$. It is obvious that the map $S \to M_g$ is not a bundle
projection, since e.g.\ for $r=1$ and $g \geq 4$ the dimension of 
$M_g$ is larger than $\dim(S)$. There is another point that deserves caution.
\par
\begin{Rem} \label{notequi}
In general, the map $p:S \to M_g$ is not equidimensional over its image.
\end{Rem}
\par
For $g \leq 3$ the map $p:S \to M_g$ is equidimensional over its image.
This is obvious for $g=2$ and the case by case discussion for $g=3$
is done in the proof of Theorem~\ref{endg3}.
\par
\begin{proof}[Proof of Remark \ref{notequi}]
An example is the stratum $S=\Omega M_5(2,2,2,2)$. This stratum has
dimension $13$ while $\dim M_5 = 12$. There is a canonical
choice of a theta characteristic $\kappa$ on $S$: If 
$$\div(\omega) = \sum_{i=1}^4 2 P_i, \quad \text{let} \quad
\kappa = \sum_{i=1}^4 P_i.$$
Over a hyperelliptic curve $X_0$ with involution $i$ one has
$h^0(\kappa) =3$, since the divisors of the form
$\cOO_{X_0}(P+i(P)+Q+i(Q))$ are linearly equivalent for all $(P,Q)$.
Since $h^0(\kappa) \mod 2$ is deformation invariant (\cite{At71}), the
fiber dimension of $p$ at the generic point of $M_5$ is one or three.
Since the latter is impossible by the dimension count, the claim
follows.
\end{proof}
\par
We let $D$ be the zero divisor of the universal holomorphic
one-form on the stratum $S$. Since a stratum $S$ is by definition disjoint 
from the hyperplane at  infinity and the zero section, $D$ may equivalently
be defined as the divisor such that the composition
\begin{equation} \label{divstrat} 
f^* \cOO_\PP(-1)|_S  \to f^* p^*(f_* \Omega^1_{X/M_g} \oplus \cOO_{M_g})|_S
\to  f^* (f_* \Omega^1_{X/\PP}|_S)
\to \Omega^1_{X/\PP}|_S \otimes \cOO_X(-D)
\end{equation}
is an isomorphism. By definition of the strata, $D$
is \'etale over its image under $f$. We denote
by $D_{\rm red}$ the corresponding reduced divisor (of
degree $r$). The bundle $\cOO_\PP(1)|_S$ is isomorphic to $\cOO_S$, 
but we keep track of it in order to naturally identify some
maps as multiplication maps (see Lemma~\ref{givenbymult}).
In the sequel we define the {\em relative tangent bundle} 
$T_{S/p(S)}$ of $p: S \to p(S)$ to be the kernel of $T_S \to p^* T_{p(S)}$ and 
similarly for restrictions of $p$ to submanifolds $B \subset S$. 
Given the Remark~\ref{notequi}, 
it would be more cautious to define it in the derived category
or to work dually with cotangent bundles. We will do the latter below.
In the applications $T_{B/p(B)}$ will always be a vector
bundle and fits well with geometric intuition. 
\par 
The following theorem is equivalent to Veech's theorem (\cite{Ve86}) on the existence
of period coordinates.
\par
\begin{Thm} \label{tanstratcoho}
The tangent bundle to a stratum can be described
by the following hypercohomology sheaf:
\begin{equation}\label{tseq}
\begin{split}
T_S &\cong R^1 f_*\left((\Omega^1_{X/S})^\vee \to
 f^* \cOO_{\PP}(1)|_S \otimes \Omega_{X/S}^1(D_{\rm red}-D)\right) \\
 &\cong \cOO_{\PP}(1)|_S \otimes R^1 f_*
\left(\cOO_S(-D) \to \Omega_{X/S}^1(D_{\rm red}-D)\right).
\end{split}
\end{equation}
In particular, 
$T_S$ carries a linear structure, which coincides with the one
defined by period coordinates in Section~\ref{affstruc}.
\par
Moreover the map 
$p:S \to p(S) \subset M_g$ induces a map
$T_S \to p^*T_{p(S)}$ with kernel 
$$ T_{S/p(S)} = \cOO_\PP(1)|_S \otimes f_* \Omega_{X/S}^1(-D+D_{\rm red}).$$
\end{Thm}
\par
\begin{proof}
Over the stratum $S$ and by definition of $D$
in equation (\ref{divstrat}) the map $\varphi$ used in 
Theorem~\ref{tangtoPP}
factors through the composition
$$ \xymatrix{\varphi_S : 
(\Omega^1_{X/S})^\vee 
\ar[r]^(.5){{\rm can}} & f^* \cOO_{\PP}(1)(-D)|_S \ar[r]^(.4){d}
&  f^* \cOO_{\PP}(1)|_S\otimes \Omega_{X/S}^1(-D+D_{\rm red}).} \\
$$
We denote by $\mfC_S^\bullet$ 
the two-term complex in degrees $0$ and $1$ given by $\varphi_S$.
As in the preceding theorem one constructs a
map $T_S \to R^1 f_* (\mfC_S^\bullet)$ as the
connecting homomorphism associated to the appropriate
exact sequence of complexes.
Over each point the stalk of $R^1 f_* (\mfC_S^\bullet)$ 
is the tangent space to $S$ by deformation theory. 
Together this implies that the two bundles are isomorphic.
\par
The second identification of $T_S$ stems from the isomorphism
of complexes 
$$ \mfC^\bullet_S \cong [f^* \cOO_{\PP}(1)|_S \otimes
(\xymatrix{\cOO_S(-D) \ar[r]^<<<<d & \Omega_{X/S}^1(D_{\rm red}-D)})]$$
induced by \eqref{divstrat} and the projection formula.
\par
The first hypercohomology of the de Rham complex 
carries a flat connection (see \cite{KaOd68}
for an algebraic definition). Furthermore the usual derivative
defines a connection on $\cOO_\PP(1)|_S$ since $S \cap \infty = \emptyset$.
The tensor product of these flat connections yields a flat connection
on $T_S$. Denote by $j: X \sms D \to X$ the inclusion. One checks
that $$\Ker(d:\cOO_S(-D) \to \Omega_{X/S}^1(D_{\rm red}-D)) = j_! \CC,$$
where $j_!$ denotes the extension by zero. Hence
$$R^1 f_*
\left(\cOO_S(-D) \to \Omega_{X/S}^1(D_{\rm red}-D)\right) \cong
 R^1 f_* j_! \CC.$$ 
Now it suffices to unwind definitions to see that the
connection just defined on $T_S$ coincides with the linear structure
by period coordinates an in Section~\ref{affstruc}.
\par
The identification of $T_{S/p(S)}$ follows from the stupid
filtration of $\varphi_S$ plus the usual Kodaira-Spencer theory
as at the end of the proof of Theorem~\ref{tangtoPP}.
\end{proof}
\par
It is easy to see the foliation by constant absolute periods in
this setting:
\par
\begin{Cor} \label{defofq1}
The inclusion $\cOO_X(-D) \hookrightarrow \cOO_X$  induces a surjection
$$ q: T_S \twoheadrightarrow \cOO_\PP(1)|_S \otimes 
R^1 f_*\left(\cOO_X  \to \Omega^1_{X/S}\right). $$
The restriction of $q$ to the relative tangent bundle $T_{S/p(S)}$ is 
injective. 
\end{Cor}
\par
This restriction
will be denoted by $q_\rel$ in the sequel.
\par
\begin{proof} The commutative diagram
\begin{equation} \label{kplxincl}
 \xymatrix{
(\Omega^1_{X/S})^\vee \ar[r]^{\rm can} \ar[d]^{\varphi_S} \,\,&  
\ar@{^{(}->}[r] f^* \cOO_{\PP}(1)(-D)|_S & f^* \cOO_{\PP}(1) \ar[d]^d\\
 f^* \cOO_{\PP}(1)|_S \otimes \Omega_{X/S}^1(-D+D_{\rm red})\,\,
\ar@{^{(}->}[rr] & \ 
& f^* \cOO_{\PP}(1)|_S \otimes \Omega_{X/S}^1}
\end{equation}
together with the projection formula defines the map $q$.
The surjectivity of $q$ can be checked pointwise at each $s \in S$.
The fiber at $s$ of the complex $\varphi_S$ is quasiisomorphic to 
$j_! \CC$, where $j: X \sms D \to X$ is the inclusion. The map
of complexes in the diagram (\ref{kplxincl}) comes from the
inclusion $j_! \CC \hookrightarrow \CC$. The spectral sequence
$$ E_2^{i,j} = H^i(X_s,H^j(\mfC^\bullet_S))$$ has non-zero
terms only for $j=0$, hence degenerates at $E_2$. 
Thus, $q$ restricts fiberwise to the map 
$$ H^1(X_s,j_! \CC) \to H^1(X_s, \CC),$$
which is obviously a surjection. 
\par
With the identifications of Theorem~\ref{tanstratcoho}, the map $q_\rel$ is 
the injection on the bottom of diagram~\eqref{kplxincl}.
\end{proof}
\par

\subsection{The dual viewpoint} \label{dualview}

We now dualize the results obtained in Section~\ref{tangstrat}.
The following diagram follows immediately from the definitions and
Serre duality. The middle and right column describe the identification
of the cotangent bundle to a stratum. The map from left column
to there is the dual of $q$ defined above. Precisely the concrete
description of $q^\vee|_{\Ker(\pi)}$, with $\pi$ defined below, 
will play a key role in the sequel. Let $\dR = [\cOO_X \to \Omega^1_{X/S}]$ 
denote the relative de Rham complex.
\begin{equation} \label{dualidentification}
\xymatrix{
\cOO_\PP(-1)|_S  \otimes \Ker(\pi) \ar[r]^{q^\vee|_{\Ker(\pi)}} 
\ar@{^{(}->}[d] & 
f_* (\Omega_{X/S}^1)^{\otimes 2} \ar[d] \ar@{->>}[r] & 
p^*\Omega^1_{p(S)}  \ar@{^{(}->}[d] \\
\cOO_\PP(-1)|_S \otimes R^1 f_* \dR \ar@{->>}[d]^{\id \otimes \pi}  \ar[r]^{q^\vee} & 
R^1 f_* ((\mfC^\vee)^\bullet_S) \ar[d]
\ar[r]^\cong &\Omega^1_S \ar@{->>}[d]\\
\cOO_\PP(-1)|_S  \otimes R^1 f_* \cOO_X(D-D_\red) 
\ar[r]^=& \cOO_\PP(-1)|_S  \otimes R^1 f_* \cOO_X(D-D_\red) 
\ar[r]^>>>>\cong & \Omega^1_{S/p(S)} \\
}
\end{equation}
The map in the complex $(\mfC^\vee)^\bullet_S = [f^* \cOO_\PP(-1)|_S \otimes 
\cOO_X(D-D_\red) \to (\Omega^1_{X/S})^{\otimes 2}]$ is the composition
of derivation and the identification \eqref{divstrat}.
The vertical arrow on top in the middle is not injective in general, 
its kernel is the image of $\cOO_\PP(-1)|_S  
\otimes f_* \cOO_X(D-D_\red)$. The map $q^\vee$ is defined by applying 
$R^1f_*(\cdot)$ to the map of complexes
\begin{equation} \label{dualkplxincl}
\xymatrix{
f^* \cOO_{\PP}(-1)\,\, \ar[d]^d  \ar@{^{(}->}[r]
& f^* \cOO_{\PP}(-1) \otimes  \cOO_X(D-D_\red) \ar[d] \\
f^* \cOO_{\PP}(-1)|_S \otimes \Omega_{X/S}^1 \ar[r]
& (\Omega^1_{X/S})^{\otimes 2} \\
}\end{equation}
induced by \eqref{divstrat}. The map $\pi: R^1 f_* \dR \to
R^1 f_* \cOO_X(D-D_\red)$, induced by the map of complexes
$$[\cOO_X \to \Omega_{X/S}^1] \to [\cOO_X(D-D_\red) \to 0],$$
completes the description of the diagram. The proof of the
projection formula implies immediately:
\par
\begin{Lemma} \label{givenbymult}
With the identifications $\cOO_{\PP}(-1)|_S \to f_*\Omega^1_{X/S}(-D)$
the map $q^\vee|_{\Ker(\pi)}$ restricted to $f_* \Omega^1_{X/S} \subset
\Ker(\pi)$ 
is the multiplication map of one-forms.
\end{Lemma}

\section{Submanifolds of strata with
linear structure: Examples} \label{linsubmf} 

In this section we describe basic properties of submanifolds of strata 
with linear structure. Recall our convention, that linear structures
are always inherited from the one on the stratum.
\par
We let $j: X \sms D \to X$ the inclusion over the stratum
of the complement of the universal divisor into the universal 
family of curves $f:X \to S$. With the preparation
made in the previous section the following criterion is
immediate. In the sequel, the subvarieties $B$ are always supposed
to be closed inside a stratum, but in general they are not compact.
\par
\begin{Thm} \label{DefofLLLemma}
If $B$ is a closed submanifold of a stratum $S$ with linear structure, 
then there are local subsystems 
$$\widetilde{\LL}_B  \subset R^1 (f|_{X_B})_* (j_! \, \CC) 
\quad \text{and} \quad \LL_B \subset R^1 (f|_{X_B})_* \CC =: \VV_B $$ 
such that $T_B \to T_S|_B$ has the image $\cOO_\PP(1)|_B \otimes
\widetilde{\LL}_B$
and the composition of $q$ and the tangent map of the inclusion
$$\psi: T_B \to T_S|_B \to \cOO_\PP(1)|_B \otimes \VV_B$$
maps to and onto $ \cOO_\PP(1)|_B \otimes \LL_B$.
\par
Conversely, suppose there is a local subsystem 
$\widetilde{\LL}_B \subset R^1 (f|_{X_B})_* (j_! \, \CC)$ 
such that $T_B \to T_S|_B$ has the 
image $\cOO_\PP(1)|_B \otimes \widetilde{\LL}_B$. 
Then $B$ is a linear manifold.
\end{Thm}
\par
\begin{proof} One direction is immediate from the definition of a 
linear submanifold. For the converse note that $T_B$ carries a flat
connection by hypothesis. This is sufficient by Lemma~\ref{critlinsubmf}.
\end{proof}
\par
\begin{Cor} \label{q1resiso}
The restriction $\psi_{rel}: T_{B/M} \to (\cOO_\PP(1)|_B 
\otimes \LL_B)$ of $\psi$ to the relative tangent bundle
is an isomorphism onto $\cOO_\PP(1)|_B  \otimes 
(\LL_B \cap f_* \omega_{X/B}(D_\red -D))$. In particular, for 
$B$ contained in the generic stratum,  the fiber dimension 
of $B \to M$ is constant.
\end{Cor}
\par
\begin{proof} Injectivity was shown in Corollary~\ref{defofq1} and 
surjectivity follows from the definition of $\LL_B$ in the preceding
theorem 
and the description of the image of the  relative tangent bundle
in Theorem~\ref{tanstratcoho}. In the generic stratum, $D_\red =D$
and the image of $\psi_{rel}$ is $(\LL_B)^{(1,0)}$. This is a vector
bundle, by Hodge theory.
\end{proof}
\par
To give a flavor of the notion we now present two familiar examples of 
submanifolds with linear structure in strata of $\Omega M_g$ in our language.
\par
{\bf Eigenforms over Teichm\"uller curves (\cite{Mo06a})}
\newline
A {\em Teichm\"uller curve} is an algebraic curve $C \to M_g$ in the
moduli space of curves, totally geodesic for the Teichm\"uller metric.
Restricting to the oriented case, they are constructed  as
$\GL^+_2(\RR)$-orbit of a pair $(X_0,\omega_0)$, called a 
{\em Veech surface}.
Recall that we abuse $M_g$ instead of replacing it by $M_g^{[n]}$.
Let $f: X \to C$ be the universal family over a Teichm\"uller curve.
By construction the bundle $f_* \omega_{X/C}$ comes with a 
distinguished  subbundle $\cLL$ such that the fibers $(X_0,\cLL_0)$ 
over any $0 \in C$ are Veech surfaces. 
Let $p:B \to C$ be the total space minus the zero section of the 
bundle $\cLL$. We hope that no confusion arises by using $p$ both
for $\Omega M_g \to M_g$ and for its restrictions.
\par
\begin{Prop}
The total space $B$ of the distinguished bundle of one-forms 
over a Teichm\"uller curve  is a closed submanifold of some 
stratum in $\Omega M_g$ with linear structure. It is defined over $\RR$.
\end{Prop}
\par
This is obvious, if we consider a Teichm\"uller curve as
the image of a $\GL_2^+(\RR)$-orbit. We reprove it, starting
from the characterization of a Teichm\"uller curve by
its variation of Hodge structures (VHS).
\par
\begin{proof} 
The VHS over Teichm\"uller curve  has a unique rank 
two summand $\LL$, defined over $\RR$, whose Higgs field is an 
isomorphism (see Section~\ref{affstruc}). The summand $\LL$
has the property that its $(1,0)$-part
$\LL^{(1,0)} \cong \cLL$ is the distinguished subbundle. We refer 
to \cite{Mo06a} for details and indicate
here only the consequences we need. We use moreover the 
fact that the lift $B$ of a Teichm\"uller curve to $\Omega M_g$
lies entirely in some stratum $S$. We are not aware of an
algebraic proof of this fact.
\par
We denote by $\LL_B$ the pullback of the local system $\LL$ on $C$ to $B$. 
We need to check that there is a local system underlying the image 
of $T_B \to T_S|_B$. 
\par
By construction of $B$ as total space of the bundle $\cLL$ we have 
$\cOO_\PP(1)|_B \cong p^* \cLL^{-1}$. The relative 
tangent bundle $T_{B/C}$ is trivialized by the section $\partial
/ \partial z$ where $z$ is a fiber coordinate.
Consider the  map $q$  defined in Proposition \ref{defofq1}.
The restriction of $q$ to the relative tangent bundle
is an inclusion
$$ T_{B/C}  \cong \cOO_B \to {\rm image} \otimes  p^*\cLL^{-1} 
\subset f_* \omega_{X/B} \otimes  p^*\cLL^{-1}. $$
We conclude that ${\rm image} \cong p^*\cLL$ as subbundles
of $f_* \omega_{X/B}$. On the other hand 
on $p^* T_C$ the map $q$ equals the composition of
$p^* T_C \to p^* \cLL^{-2} \subset R^1 f_* \omega_{X/B}^{-1}$ 
and the restriction of the natural map
$$R^1 f_* \omega_{X/B}^{-1} \to R^1 f_* \cOO \otimes \cOO_\PP(1)|_B
\cong R^1 f_* \cOO \otimes p^* \cLL^{-1}. $$
Both composition factors are maps that arise as $p^*$ of maps between 
bundles on $C$: For the second map this is obvious from the definition
and the first map is the Kodaira-Spencer map that depends only on
the deformation of the underlying curve, not on the additional one-form.
Moreover, all bundles and maps extend  over the compactification
$\ol{C}$  and yield a composition map
$$ T_{\ol{C}} \to R^1 f_* \omega_{X/C}^{-1} \to  R^1 f_* \cOO \otimes 
\cLL^{-1}. $$
The property 'maximal Higgs' implies that this composition is injective
onto $\cLL^{-1} \otimes \cLL^{-1}$ and splits. 
\par
Together we have shown that $q$ restricts to an isomorphism
$$\psi: T_B  \to \LL_B  \otimes \cOO_\PP(1)|_Y \subset  \VV_B 
\otimes \cOO_\PP(1)|_Y. $$
Taking an unramified cover, we may assume that the zero divisor $D \to B$
consists of sections $s_i$ of $f$. By \cite{Mo06b} Corollary~3.4
the difference $s_i - s_j$ of any two sections is torsion in 
${\rm Pic}^0(X/B)$.
This implies that the relative periods are rational multiples of
the absolute periods. Equivalently, the local
system $\LL_B$ lifts uniquely to a local system 
$\widetilde{\LL}_B$
such that $T_B \to {\LL}_B\otimes \cOO_\PP(1)|_B $ factors through 
$T_B \to \widetilde{\LL}_B\otimes \cOO_\PP(1)|_B$.
Obviously $\widetilde{\LL}_B$ is again defined over $\RR$ and satisfies
what is needed to apply Theorem~\ref{DefofLLLemma}. 
\end{proof}
\par
{\bf Covering constructions}
\par
It is well-known that from a $\GL^+_2(\RR)$-invariant manifolds
one can construct a new $\GL^+_2(\RR)$-invariant manifold in higher
genus by a Hurwitz space construction. The dimension will increase
by the number of ramification points. We reprove this fact, 
showing that it does not depend on the linear structure being
defined over $\RR$.
\par
\begin{Prop}
Let $B \subset S:=\Omega M_g(k_1,\ldots,k_r)$ be a manifold with linear
structure and
$f: X \to p(B)$ be the universal family over $p(B)$. Let $H \to p(B)$ be
the Hurwitz space parameterizing covers of fibers of $f$ ramified over 
$s$ points and of some fixed type. Then $B' := H \times_{p(B)} B$ is a linear
submanifold of some stratum $S':= \Omega M_{g'}(k'_1,\ldots,k'_{r'})$.
\end{Prop}
\par
\begin{proof}[Sketch of proof:]
Let  $f: X' \to B'$ be the pullback of the universal family,
let $\pi: Y' \to X'$ be the universal covering of $X'$ and
let $g = f \circ \pi$. We let also $j_Y: Y' \sms D' \to Y'$
be the inclusion of the complement of the universal divisor.
\par
By hypothesis there is a linear subsystem 
$\widetilde{\LL}_B \subset R^1 f_* (j_X)_! \,\CC$ such that
$$T_B \cong \cOO_\PP(1)|_B \otimes \widetilde{\LL}_B \subset T_S|^{\ }_B.$$
\par
We let $\LL_B$ be the image of $\widetilde{\LL}_B$ in $R^1 f_* (j_X)_! \,\CC$.
We have an identification 
$$T_{S'}|^{\ }_{B'} \cong  \cOO_\PP(1)|_{B'} \otimes R^1 g_* (j_Y)_!\, \CC
\cong \cOO_\PP(1)|_{B'} \otimes R^1 f_* ((j_X)_! \, \pi_* \CC).$$
Via the split inclusion $\CC \to \pi_* \CC$ we identify 
$ R^1 f_* (j_X)_! \, \CC \subset  R^1 f_* ((j_X)_! \, \pi_* \CC)$
as a direct summand. Moreover, we write $D=\sum_{i=1}^r k_i D_i$ and
$D'= \sum_{i=1}^{r'} k'_i D'_i$ for the universal divisor and we use
$\CC_{k_i D_i}$ to denote the skyscraper sheaf
of length $k_i$ along $D_i$. For the same reason as above, the inclusion
$f_*(\oplus_{i=1}^r \CC_{k_i D_i}) \to f_*(\oplus_{i=1}^{r'} \pi_* \CC_{k'_i D'_i})$
is split and we let $\MM$ denote the complement.
\par
It is now easy to check that 
$$T_{B'} \cong \cOO_\PP(1)|_{B'} \otimes (\widetilde{\LL}_{B'} \oplus \MM)
\subset T_{S'}|^{\ }_{B'}. $$
This shows that $B'$ has a linear structure. 
\end{proof}
\par

\section{Linear structures and 
endomorphisms of the Jacobian} \label{secHypEnd}

The aim of this section is to show one of the central results: 
\par
\begin{Thm} \label{comphypend}
Suppose that $B$ is a closed algebraic submanifold of the generic
stratum $S=\Omega M_g(1,\ldots,1)$ with linear structure. Then
there are three possibilities:
\begin{itemize}
\item[(i)] $B$ is a connected component of $S$.
\item[(ii)]  $g \geq 3$ and $B$ is the preimage in $S$
of the hyperelliptic locus in $M_g$. 
\item[(iii)] $B$ parameterizes curves with a Jacobian whose endomorphism ring
is strictly larger than $\ZZ$.
\end{itemize}
\end{Thm}
\par
Case iii) also contains the covering constructions, see the end
of the previous section, since the Jacobian splits in these cases.
The non-generic strata need closer inspection. We do this case
by case in genus $g=3$.
\par
In this section we will rely on Deligne's semisimplicity for
the VHS over the manifolds under consideration. That is why
we need $B$ to be algebraic or, a priori
slightly weaker, that any plurisubharmonic function on $B$ which is bounded
above is in fact constant. 
\par 
We continue with the notation used in the previous section, in 
particular in Theorem~\ref{DefofLLLemma}.
\par
\begin{Lemma} \label{notallthenEnd}
If $\LL_B \varsubsetneqq \VV_B$ then $B$ parameterizes curves 
with a Jacobian whose 
endomorphism ring is strictly larger than $\ZZ$.
\end{Lemma}
\par
\begin{proof} By semisimplicity (\cite{De87} Proposition 1.13) the 
$\CC$-VHS on the local system $\VV_B$  decomposes as
\begin{equation} \label{VHSdecomp}
\VV_B = \bigoplus_{i \in I} (\LL_i \otimes W_i) 
\end{equation}
with irreducible local systems $\LL_i$ of rank $d_i$
carrying a weight one VHS, whose Higgs field is non-zero,
and with $W_i$ vector spaces carrying a weight zero VHS.
\par
We claim that, as in  \cite{ViZu04}  Lemma 3.2, each of the 
$\LL_i$ may be chosen to be defined over some number field $L_i$.
In fact, let ${\mathcal G}(d_i,\VV_B)$ be the set of local systems in $\VV_B$
of rank $d_i$ and let ${\rm Grass}(d_i,V)$ be the Grassmann
variety of subspaces of dimension $d_i$ of a the vector space 
$V$, which is defined as the fiber of $\VV_B$ at some base point.
The set ${\mathcal G}(d_i,\VV_B)$
is naturally a subvariety of 
$${\rm Grass}(d_i,V) \times_{{\rm Spec}\, \QQ} {\rm Spec}\, \CC$$
consisting of $\pi_1(B)$-invariant points.
The subset in ${\mathcal G}(d_i,\VV_B)$ of rank $d_i$ local subsystems,
whose projection to $\LL_i$ along its complement under the
polarization is non-zero, is Zariski-open. It hence contains
an element, which is defined over a finite extension $L_i$ of $\QQ$. 
By the irreducibility of $\LL_i$ the non-zero map must in fact
be an isomorphism. 
\par
From the hypothesis $\LL_B \varsubsetneqq \VV_B$ we deduce that
the decomposition (\ref{VHSdecomp}) is non-trivial. Let $K$
be the Galois closure of $L_1/\QQ$. 
Then for all $\sigma \in \Gal(K/\QQ)$
the local systems $\LL_1^\sigma$ are also among the $\LL_i$. 
Let $L\subset K$ be the
field fixed by all automorphisms $\sigma$ such that
$\LL_1^\sigma \cong \LL_1$ as sub-local systems of $\VV$.  
For $a \in L$ the map 
$$\bigoplus_{\sigma \in \Gal(K/\QQ)/\Gal(K/L)}  
\sigma(a) \cdot{\rm id}_{\LL_{i(\sigma)} \otimes W_{i(\sigma)}} 
\in \End(\VV_B),
\quad \text{where}\,\, \LL_{i(\sigma)}:=\sigma(\LL_1) 
$$
This endomorphism is of bidegree $(0,0)$ and defined over $\QQ$, 
hence a non-trivial $\QQ$-endomorphism of the family of Jacobians over $B$. 
If $L = \QQ$ and $|I|>1$ then $\QQ \cdot {\rm id}_{\LL_1 \otimes W_1}$
are non-trivial $\QQ$-endomorphisms. Finally, if  $L = \QQ$ and $|I|=1$
then $W_1$ has dimension at least two. Choose a one-dimensional 
subvectorspace $W_0 \subset W_1$ and take the non-trivial endomorphisms
$\QQ \cdot {\rm id}_{\LL_1 \otimes W_0}$.
\end{proof}
\par
\begin{proof}[Proof of Theorem \ref{comphypend}]
Suppose that neither (i) nor (iii) holds for $B$.
By Corollary~\ref{q1resiso} $T_B \to (p^*T_{M_g})|_B$
is a locally split surjection of vector bundles. Hence
$(p^* \Omega_{M_g}^1)|_B \to \Omega^1_B$ is a locally split inclusion.
We can say more here:
Following the notations of Theorem~\ref{DefofLLLemma} and by 
Lemma~\ref{notallthenEnd}
the map 
$$ \psi: T_B \to \cOO_\PP(1)|_B \otimes \VV_B$$
is surjective. A first consequence of this fact is that over each 
point in $p(B)$ the fibers of $p|_S$ and $p|_B$ coincide. Thus, there is 
locally over any sufficiently small open set
$U \subset M_g$, a splitting $s_U:\Omega^1_{p(B)}|_U \to 
\Omega_{M_g}^1|_U $ of
the surjection $\Omega_{M_g}^1|_{p(B)} \to  \Omega^1_{p(B)}$.
\par
Second, the composition
$$p^* T_B \to p^*T_{M_g} |_B = p^*R^1 f_* (\Omega^1_{X/M_g}|_B)^{-1}
\to \cOO_\PP(1)|_B \otimes R^1 f_* \cOO_{X_B}$$
has to be surjective, too. Using the notation of Lemma~\ref{givenbymult},
we have $f_* \Omega^1_{X/S} = \Ker(\pi)$, since all zeroes are simple on
the generic stratum.
Dualizing this composition yields that the multiplication map
$$ \psi^\vee|_{\Ker(\pi)}: \cOO_\PP(-1)|_B \otimes (p^* f_* \Omega^1_{X/M_g})|_B
\to (p^*f_* (\Omega^1_{X/M_g})^{\otimes 2})|_B \to p^*\Omega^1_{p(B)}$$ has to be 
injective and factors, over say $p^{-1}(U)$, 
through $\cFF_U:=p^*(s_U(\Omega^1_{p(B)}))$.
By Lemma~\ref{givenbymult} this implies that for all line bundles
$\cEE_U \subset f_* \omega_{X/M_g}|_U$ generated by a differential
with simple zeroes the image of the multiplication map
$$\cEE_U \otimes f_* \omega_{X/M_g}|_U \to f_* \omega_{X/M_g}^{2}|_U$$
lies in the subbundle $\cFF:=p^*(s_U(\Omega^1_{p(B)}))$. Since
one-forms with simple zeroes span the space of holomorphic one-forms, 
this implies that the whole image of the multiplication map
$$f_* \omega_{X/M_g}|_U \otimes f_* \omega_{X/M_g}|_U \to 
f_* \omega_{X/M_g}^{2}|_U$$
lies in $\cFF$.
\par 
By M.Noether's theorem (e.g.\ \cite{ACGH85} III.\S~2), for $g=2$
or at a point outside the hyperelliptic locus, the multiplication
map is surjective. Hence if $p(B)$ is not contained in the
hyperelliptic locus, $\cFF = p^*(s_U(\Omega^1_{p(B)}))$ and
$B$ coincides with the whole stratum.
\par
If $p(B)$ is contained in the hyperelliptic locus, the 
image of the multiplication map is  a subbundle in 
$f_* \omega_{X/M_g}^{\otimes 2}$ of rank $2g-1$ as can
be checked easily on an explicit basis of one-form. See
Theorem~\ref{endg3} for more involved arguments in the same style.
For dimension reasons, $p(B)$ has to coincide with the
hyperelliptic locus.
\par 
Finally we remark that this case actually occurs, i.e.\ that the
bundle $B$ over the hyperelliptic locus has a linear structure
defined over $\RR$. For this purpose we check that 
$B$ is $\GL_2^{+}(\RR)$-invariant. This easily follows from the
fact that $-1$ is in the center of $\GL_2^{+}(\RR)$ or, alternatively,
that $B$ is the image of stratum of quadratic differentials
with $2g+2$ simple poles and $g-1$ double zeros on $\PP^1$
under a canonical double covering map, see e.g.\ \cite{KoZo03}~Lemma 1. 
\end{proof}
\par
\par
\begin{Thm} \label{endg3}
For a closed submanifold $B$ with linear structure in a stratum
$S \neq \Omega M_3(1,\ldots,1)$ of $\Omega M_3$ there are the 
following possibilities:
\begin{itemize}
\item[i)] $B$ is a connected component of $S$.
\item[ii)] $B$ parameterizes curves whose Jacobian has an endomorphism
ring strictly larger than $\ZZ$.
\item[iii)]  $B$ is the preimage of the hyperelliptic locus in
$\Omega M_3(2,1,1)$ or in  the unique component $\Omega M_3(2,2)^{\rm odd}$
of $\Omega M_3(2,2)$
 that does not exclusively consist of hyperelliptic curves.
\item[iv)] $B$ has dimension $5$, is disjoint from the
hyperelliptic locus and $B \subset \Omega M_3(2,2)^{\rm odd}$
or $B \subset \Omega M_3(3,1)$.
\end{itemize} 
\end{Thm}
\par
We do not claim, that the manifolds in iv) exist. In fact
we doubt their existence, but cannot not prove it.
\par
\begin{proof}
{\it Case $S=\Omega M_3(4)$:} This stratum has two connected components
one contained in the hyperelliptic locus, the other one 
disjoint from the hyperelliptic locus. In both cases the fibers of $p|_S$
are one-dimensional, since zeros of order $\geq 3$ on a curve of genus $3$
are Weierstra\ss\ points, hence there are only finitely many of them
on a fixed curve. Consequently, the image of both connected 
components in $M_3$ is
$5$-dimensional. If $B\subset S$ belongs not to case ii) then
$q$ is surjective by Lemma \ref{notallthenEnd}. In particular the 
induced map on $p^* T_{p(S)}|_B$
surjects onto the whole of
$(\cOO_\PP(1) \otimes R^1 f_*(\cOO_B \to \omega^1_{X/B}))/q_\rel(T_{S/p(S)})$.
This bundle has rank $5$. Hence $\dim p(B)=5$ and $B$ belongs to case i).
\par
{\it Case $S=\Omega M_3(3,1)$:} This stratum does not intersect
the hyperelliptic locus since the hyperelliptic involution would have
to fix the zeros. But on a hyperelliptic curve differentials have zeros
at Weierstra\ss\ points of even order only. Hence the canonical model
of a curve in $S$ is a plane quartic. The fibers of $p$ are 
one-dimensional: A zero of order three corresponds a Weierstra\ss\ point,
in fact to an inflexion
line of the plane quartic. The simple zero is the 4th point of 
intersection of this line with the quartic, hence fixed once the
Weierstra\ss\ point is chosen. Since $\dim(S)=7$ we deduce that
$p(S)$ is of dimension $6$, dense in $M_3$.
\par
Suppose that $B \subset S$ does not belong to case ii). 
This implies that 
$$\ol{\psi}: p^*T_{p(B)} \to (\cOO_\PP(1) \otimes R^1 
f_*(\cOO_B \to \omega^1_{X/B}))/q_\rel(T_{S/p(S)})$$
is surjective and that $B$ has dimension at least $5$. Hence $B$
belongs to case i) or iv).
\par
{\it Case $S=\Omega M_3(2,2)$:} By \cite{KoZo03} Theorem 2, there 
are two components of this stratum, both of dimension $\dim S =7$:
the component $S^{\rm hyp}:=\Omega M_3(2,2)^{\rm hyp}$ consisting 
of hyperelliptic curves
such that the hyperelliptic involution interchanges the two zeroes,
and $S^{\rm odd}:=\Omega M_3(2,2)^{\rm odd}$, the component with 
odd spin structure (see \cite{KoZo03} for the definition). 
\par
Consider $B \subset S^{\rm hyp}$ first. Since the hyperelliptic locus 
has dimension $5$ and since the generic differential on any curve
is not of type $(2,2)$
all fibers of $S^{\rm hyp} \to p(S^{\rm hyp})$ are two-dimensional
and the image  $p(S^{\rm hyp})$ is dense in the hyperelliptic locus. If
$B$ does not belong to case ii) we use the same argument as 
in Theorem~\ref{comphypend} to show that $B=S^{\rm hyp}$: 
We have to show that the images of multiplication map
$$\cEE_U \otimes f_* \omega_{X/M_g}|_U \to f_* \omega_{X/M_g}^{2}|_U,$$
for all line bundles $\cEE_U$ generated by a differential with 
two double zeroes that
are interchanged by the hyperelliptic involution, generate
the subspace of $f_* \omega_{X/S}^2$ acted on by $+1$ under the
hyperelliptic involution. On a hyperelliptic genus $3$ curve
$$ y^2 = \prod_{i=1}^7 (x-x_i) \quad \text{with}\, x_i \neq x_j 
\quad \text{for} \quad i \neq j$$
a basis of holomorphic one-forms is $x^idx/y$ for $i=0,1,2$.
The $+1$-eigenspace is generated by $x^i(dx/y)^2$ for $i=0,\ldots,5$. 
The one-forms
$(x-x_0)^2dx/y$ for $x_0 \neq x_i$, $i=0,\ldots,7$ all have the
prescribed type of zeros. Obviously the vector space they generate
consists of all  holomorphic one-forms. 
Hence their products
with the  $x^idx/y$ for $i=0,1,2$ indeed generates the $+1$-eigenspace. 
\par
Consider now $B\subset S^{\rm odd}$. A point in this stratum
defines an odd theta characteristic, as in the proof of Remark~\ref{notequi}.
There are only finitely many, in fact 28, such theta characteristics
on each curve. Moreover the space of global sections of a  fixed
theta characteristic cannot be three-dimensional. Hence it is
one-dimension for all curves in $p(S^{\rm odd})$ by \cite{At71}.
It follows that the fibers of $p:S^{\rm odd} \to p(S)$ are all 
one-dimensional. 
\par 
As above, if
$B$ does not belong to case ii),  the map $\ol{\psi}$ induced by 
$\psi$ on $p^* T_{p(B)}|_B$ surjects onto the whole of
$(\cOO_\PP(1) \otimes R^1 f_*(\cOO_B \to \omega^1_{X/B}))/q_\rel(T_{S/p(S)})$.
This bundle has now again rank~$5$. Hence $B$ belongs to case
i) or iv).
\par
{\it Case $S=\Omega M_3(2,1,1)$:} This stratum has dimension $8$
and maps surjectively to $M_3$ with $2$-dimensional fibers.
In order to argue as in the proof of Theorem~\ref{comphypend}
we have to show that for a genus $3$ curve the
one-forms of type $(2,1,1)$ generate the space of holomorphic
one-forms. For hyperelliptic curves this is follows as in the
case $\Omega M_3(2,2)^{{\rm odd}}$. For a smooth plane quartic, any
tangent line except for finitely many bitangents and inflexion
lines corresponds to a one-form of type $(2,1,1)$. The tangent
lines to a quartic form a (singular) curve $C^*$ of degree $12$ in 
the dual of $\PP^2$. If the generation result we need was false, 
$C^*$ was contained in, hence equal to, a hyperplane in $(\PP^2)^\vee$.
This is absurd.
\end{proof}
\par

\section{Hilbert modular varieties}

Fix a totally real number field $K$ with $[K:\QQ]=g$ and
an order $\fraco \subset K$. We denote by $\fraco^\vee$
the inverse different and we let
$$ \SL(\fraco \oplus \fraco^\vee)= \left\{ \left(
\begin{matrix} a & b \\ c & d\\ \end{matrix} \right) \in \SL_2(K): \,
a \in \fraco, \,b\in (\fraco^\vee)^{-1}, \,c \in \fraco^\vee,\,
d \in \fraco \right\}.$$ 
The Hilbert modular variety $\HH^g/\SL(\fraco \oplus \fraco^\vee)$
parameterizes principally polarized abelian $g$-folds $A$ together 
with the choice 
of some real multiplication, i.e.\ with a map $\fraco \to \End(A)$.
Forgetting this choice yields a map from 
$\HH^g/\SL(\fraco \oplus \fraco^\vee)$ to the moduli space 
of principally polarized abelian $g$-folds. This map ramifies
along diagonals in $\HH^g$ unless we keep track of the stack
structure of $A_g$. We avoid this by passing right away to
a suitable level structure:
\par
Principally polarized abelian $g$-folds $A$ together 
with the choice of real multiplication and a choice of a 
basis of $\fraco \oplus \fraco^\vee$, symplectic for
the pairing
\begin{equation}\label{tracepairing}
\langle(x,y),(x',y')\rangle := {\rm tr}(xy'-x'y),
\end{equation}
are parameterized by $\HH^g$. Forgetting the choice of real multiplication
and fixing a $\ZZ$-basis $a_1,\ldots,a_g$ of $\fraco$ yields
an embedding
$$j: \HH^g \to \HH_g, \quad (\tau_1,\ldots,\tau_g) \mapsto
M^{T} \cdot \diag(\tau_1,\ldots,\tau_g) \cdot M,$$
where $\HH_g$ denotes the Siegel upper half plane, where $M = (\sigma_j(a_i))$ 
and $\sigma_j:K \to \RR$ are the
$g$ real embeddings of $K$. The map $j$ is equivariant with
respect to the action of $\SL(\fraco \oplus \fraco^\vee)$ on 
$\HH^g$ and the action on $\HH_g$ of its image in $\Sp(2g,\ZZ)$, 
obtained by expressing the entries of 
$\SL(\fraco \oplus \fraco^\vee)$ with respect to $a_1,\ldots,a_g$
and its dual basis with respect to the pairing \eqref{tracepairing}.
\par
Let $\Sp^{[n]}(2g,\ZZ):=\Ker(\Sp(2g,\ZZ) \to \Sp(2g,\ZZ/n\ZZ))$.
For $n \geq 3$ the quotient map $\HH_g \to \HH_g/\Sp^{[n]}(2g,\ZZ)$ is unramified
and if we let 
$$\SL^{[n]}(\fraco \oplus \fraco^\vee) = \SL(\fraco \oplus \fraco^\vee)
\cap \Sp^{[n]}(2g,\ZZ)$$
and $A_g^{[n]} := \HH_g/\Sp^{[n]}(2g,\ZZ)$ 
we obtain the two right columns of the following diagram.
\begin{equation}\label{levelstructure}
\xymatrix{
M_g^{[n]} \ar@/^/[rrrdd]^<<<<<{t^{[n]}} \ar[d]
& \HH^g \ar[rr]^{j} \ar[dd]&& \HH_g \ar[dd] \\
M_g^{[n]}/\pm 1 \ar@/^/[rrrd]^<<<<<{t^{[n]}/\pm 1} \ar[dd]\\
 &
Y:=\HH^g/\SL^{[n]}(\fraco \oplus \fraco^\vee)  \ar[rr]^{j^{[n]}} \ar[d]
&& A_g^{[n]}  \ar[d]\\
M_g \ar@/_/[rrr]& \HH^g/\SL(\fraco \oplus \fraco^\vee) \ar[rr] && A_g \\
} 
\end{equation}
In particular the map $j^{[n]}$ is an embedding. We will from
now on call, slightly abusively, $Y$ a {\em Hilbert modular variety}.
\par
Let $\ol{Y}$ be a smooth compactification with $S:=\ol{Y}\sms Y$
a normal crossing divisor. We specify the compactification  below.
Recall from \cite{ViZu07} the decomposition
of the VHS over $Y$. The logarithmic cotangent bundle decomposes as
$$\Omega^1_Y (\log S) 
\cong \bigoplus_{i=1}^g \Omega_i $$
defined as follows. The restriction of $\Omega_i$ to $Y$ is generated
by $dz_i$ where the $z_i$ are coordinates of factors of the universal
covering of $Y$. If the compactification is chosen as in \cite{Mu77},
the bundles $\Omega_i|_Y$ extend to $\ol{Y}$ and their sum
is the logarithmic cotangent bundle (Theorem~3.1 and Proposition~3.4
in loc.~cit.).
\par
Let $g: A \to Y$
be the universal family over $Y$. Then the VHS decomposes as
\begin{equation} \label{HMVHS}
\VV_\CC := R^1g_* \CC = \bigoplus_{i=1}^g \LL_i, 
\end{equation}
where the $\LL_i$ are rank two local systems, in fact defined over $K$
and Galois conjugate. Moreover all the $\LL_i$ are maximal Higgs, 
that is the Higgs field
$$\theta: \VV_\CC^{(1,0)} \to  \VV_\CC^{(0,1)} \otimes \Omega^1_Y (\log S)$$
is the direct sum of isomorphisms
\begin{equation} \label{HilbVHS}
\theta_i: \cLL_i:=\LL_i^{(1,0)} \to \LL_i^{(0,1)} \otimes \Omega_i
\end{equation}
composed with the natural inclusions $\Omega_i \to \Omega^1_Y 
(\log S)$. Recall also that $\LL_i^{(0,1)} = \cLL_i^{-1}$ and
hence $\Omega_i \cong \cLL_i^{\otimes 2}$. 
\par
From this can describe the tangent map to the inclusion $t^{[n]}$.
Choose a smooth compactification $\ol{A_g}$ of $A_g$ with boundary $S'$ and
choose $\ol{Y}$ to be the closure of $Y$ in $\ol{A_g}$.  
\par
\begin{Lemma} \label{KSHM}
The tangent map to $t^{[n]}$
$$ T_Y (-\log S) \cong  \bigoplus_{i=1}^g \Omega_i^\vee
\to T_A(-\log S')|_Y  \cong (\Sym^2(\bigoplus_{i=1}^g \cLL_i^{-1}))^\vee $$
is given by the injection onto the summands $\oplus_{i=1}^g \cLL_i^{-2}$.
More precisely, there is no non-zero map 
$$m_{i,j,k}: \cLL_i \otimes \cLL_j \to \cLL_k^2$$
unless $k=i=j$, in which case this map is an isomorphism.
\end{Lemma}
\par
\begin{proof}
The first statement follows from the second and the second 
follows from \cite{ViZu07} Lemma 1.6: The
$\cLL_i$ are line bundles with the same slope with respect to
$\omega_Y(\log S)$. Hence maps of the form $m_{i,j,k}$ are either 
zero or an isomorphism. Comparing $c_1(\cLL_i)$ we deduce from loc.~cit.
that the latter cannot happen. 
\end{proof}
\par
\subsection{Linearity of Hilbert modular surfaces: Reproving \cite{McM03a}}

Let $Y_M$ be the preimage of a Hilbert modular surface
in $M_2^{[n]}$. Since all curves with $g=2$ are hyperelliptic, 
this is an unramified double covering of $Y$ and all statements
of Higgs fields or Kodaira-Spencer maps being isomorphisms remain true 
on $Y_M$. Let $p:B^{\rm total} \to Y_M$ be the {\em eigenform locus} over 
$Y_M$, i.e.\
the subset of $\Omega M_2$ consisting of pairs $(X,\omega)$ 
where  $\fraco$ acts as $a \cdot\omega = \iota(a)\omega$ 
for a fixed embedding $\iota:K \to \RR$ for all $a \in \fraco$.
\par
The aim of this section is to reprove Theorem~{7.1} in \cite{McM03a}
stating that $B^{\rm total}$ is $\GL_2^+(\RR)$-invariant. In fact
we will show:
\par
\begin{Prop} \label{HM2linear}
The eigenform locus $B := B^{\rm total} \cap \Omega M_2(1,1)$ of a
Hilbert modular surface in the stratum 
$\Omega M_2(1,1)$ is a closed submanifold with linear structure
defined over $\RR$.
\end{Prop}
\par
The full strength of Theorem~{7.1} in loc.\ cit.\ follows from
Proposition~\ref{Kontconverse} and a local consideration around
$B^{\rm total} \cap \Omega M_2(2)$, the Teichm\"uller curves
discovered in \cite{McM03a} Theorem~1.3.
\par
\begin{proof} 
We need to show that there is a rank two local system 
$\LL_B \subset \VV_B$ defined over $\RR$
such that $\psi(T_B) = \cOO_\PP(1)|_B \otimes \LL_B$.
Then, since the fibers of $q$ are one-dimensional on $\Omega M_2(1,1)$ 
and since $B$ is three-dimensional, the image of $T_B$ in $T_S|_B$
coincides with the full preimage of $\cOO_\PP(1)|_B \otimes \LL_B$. 
It thus carries a flat structure. The conclusion then follows
from Theorem~\ref{DefofLLLemma}.
\par
Changing the enumeration of the summands if necessary,  we may suppose that
$a \in \fraco$ acts on the VHS via $a \cdot {\rm id}_{\LL_1}
\oplus \sigma(a) \cdot {\rm id}_{\LL_2}$, where $\sigma$ generates
$\Gal(K/\QQ)$. Hence by definition $\cOO_\PP(1)|_B \cong p^* \cLL_1^{-1}|_B$
and $\psi$ restricts to
$$ \xymatrix {\psi_\rel: T_{B/Y_M} \cong \cOO_B \ar[r]^<<<<\sim 
& \cLL_1|_B \otimes \cOO_\PP(1)|_B \subset
f_* \omega_{X/B} \otimes \cOO_\PP(1)|_B}.$$
\par
It remains to show that the quotient 
$$\ol{\psi}: p^*T_{Y_M} \to R^1 f_* \omega_{X/B}^{-1} \to
R^1 f_* \cOO_X \otimes \cOO_\PP(1)|_B$$
has an image equal to 
$p^*\cLL_1^{-1}|_B \otimes \cOO_\PP(1)|_B$. All bundles and
maps involved in the composition arise as pullback $p^*$ of
bundles and maps on $p(B) \subset Y_M$. In fact the first map is
the Kodaira-Spencer map and the second is the dual of the
multiplication map. Hence we need to show that
$$ T_{p(B)} \to R^1 f_* \omega_{X/p(B)}^{-1} \to
R^1 f_* \cOO_X \otimes \cLL_i^{-1}|_{p(B)}$$
maps onto $ \cLL_i^{-1}|_{p(B)}\otimes \cLL_i^{-1}|_{p(B)}$. Dualizing 
yields the map
$$\cLL_i|_{p(B)} \otimes f_* \omega_{X/p(B)} \to 
{\rm Sym}^2 (f_* \omega_{X/p(B)}) = \Omega^1_{A_2^{[n]}}|_{p(B)}
\to f_* \omega_{X/p(B)}^{\otimes 2} = \Omega^1_{M_2^{[n]}}|_{p(B)}
\to \Omega^1_{p(B)}$$
which is the multiplication of one-forms composed with the 
dual Kodaira-Spencer map. Ignoring the factorization through
$\Omega^1_{M_2^{[n]}}|_{p(B)}$ we know this map from Lemma~\ref{KSHM}:
It is the natural inclusion composed with the projection along 
the direct summand 
$$\oplus_{i \neq j} \cLL_i \otimes \cLL_j \subset
{\rm Sym}^2 (f_* \omega_{X/p(B)}).$$ 
\par
Dualizing again we conclude that $\psi$ is as desired.
\end{proof}
\par
\subsection{Obstacles to linearity for 
Hilbert modular threefolds}

In this section we analyze to which extent the linearity results
for Hilbert modular surfaces apply to higher dimensions. For $g \geq 4$
the image of $M_g$ is no longer dense in $A_g$. It fact, it is
shown in \cite{dJZh06} that Hilbert modular $g$-folds lie generically
outside $M_g$ for $g \geq 4$ with possible exceptions for $g=4$.
For $g=3$ the situation is better:
\par
\begin{Lemma} \label{redcodim2}
A Hilbert modular threefold $Y \subset A_3^{[n]}$ lies generically
in $t^{[n]}(M_3^{[n]})$. The complement has codimension
at least two.
\end{Lemma}
\par
\begin{proof}
For $g=3$ the image of $t^{[M_3]}$ is dense in $A_3$, the complement
$W$ of the abelian threefolds, that are reducible as principally
polarized abelian varieties. Suppose that $A \cong A_1 \times A_2$ is a 
generic fiber in a component $W_0$ of $W \cap Y$.
We number the factors such that $\dim A_i = i$. 
If $A_2$ was a simple abelian variety or if
none of factors of $A_2$ is isogenous to $A_1$, then 
$$\End_\QQ(A) = \End_\QQ(A_1) \times \End_\QQ(A_2) $$
and $\End_\QQ(A_2)$ is a $\QQ$-algebra of rank $2$ or $4$.
The same argument applies if $A_2$ is isogenous to $E_1 \times E_2$ and
say $E_1$ is isogenous to $A_1$ while $E_2$ is not. 
Hence we may assume $A_2$ is  isogenous to $E_1 \times E_2$
where $E_1$, $E_2$ and $A_1$ are all isogenous. But this limits
the possibilities for $A_2$ given $A_1$ to a countable number
and we conclude that $\dim W_0 \leq 1$.
\end{proof}
\par
As in the case $g=2$, let $Y_M$ be the preimage of $Y$ in
$M_3^{[n]}$ and let $\Omega M^{[n]}_g \supset B^{\rm total} \to Y_M$ 
be the $\CC^*$-bundle of eigenfoms with respect to a fixed
embedding of $K \to \RR$. Let $\Omega M_3(\fraco)$ be the
stratum the generic eigenform maps to and denote by 
$B$ the intersection of $B^{\rm total}$ and $\Omega M_3(\fraco)$.
\par
\begin{Q}
What is this stratum $\Omega M_3(\fraco)$, does it depend on $\fraco$ ?
\end{Q}
\par
For $\fraco=\ZZ[\zeta_7+\zeta_7^{-1}]$ the stratum $\Omega M_3(\fraco)$ 
is the generic stratum $\Omega M_3(1,1,1,1)$. This can be checked
using the explicit equation given in \cite{Lo05} of the
Teichm\"uller curve in Veech's series (\cite{Ve89})
generated by $y^2=x^7-1$ and the one-form $dx/y$.
\par
\begin{Thm} \label{obstacles}
Let  $B:=B_\fraco$ denote the locus of eigenforms $B:=B_\fraco$ 
over a Hilbert modular threefold and let $\Omega M_3(\fraco)$ be the
stratum of $\Omega M_3$ the generic eigenform maps to.  
The manifold $B$ inherits from $\Omega M_3(\fraco)$ a linear structure 
if and only if
\begin{itemize}
\item[i)] there  is {\em no} point $b \in B$ corresponding to a
hyperelliptic curve, such that the normal bundle to the hyperelliptic
locus at $b$ coincides (as subbundle of the cotangent bundle) with one of the 
distinguished subbundles $\Omega_i$, the cotangent directions to 
the natural foliations of a Hilbert modular threefold {\em and}
\item[ii)] either $\Omega M_3(\fraco) = \Omega M_3(2,1,1)$ or
$\Omega M_3(\fraco) = \Omega M_3(1,1,1,1)=:S$ and the local system
defined by the zeros of the eigenform extends the pullback of
$\LL_1$ to $B$ to a local subsystem of $R^1f_* (j_! \CC)|_{Y_M}$.
\end{itemize}
\end{Thm}
\par
Before proving this, several remarks are necessary. First,
concerning the first condition we have:
\par
\begin{Prop} \label{hypintHM3}
The intersection of $Y_M$ with the hyperelliptic locus $H$ in
$M_3^{[n]}$ is non-empty, $Y_M$ is not contained in $H$ 
and some components of 
$H \cap Y_M$ are codimension one in $Y_M$.
\end{Prop}
\par
\begin{proof}
Take the Baily-Borel-Satake compactification of $Y$. This
gives a projective embedding such that the boundary consists 
of codimension $\geq 2$ components. Hence an intersection
with general hyperplanes produces a curve $C_Y \subset Y$
that avoids both the boundary and, by Lemma~\ref{redcodim2}, 
the reducible locus.
Consequently, $C:=(t^{[n]})^{-1}(C_Y)$ is a compact curve lying
entirely in $Y_M$. 
\par
Suppose that the intersection $H \cap Y_M$ is empty. The
class of $H \in {\rm Pic}(\ol{M_3})$ is a positive multiple
of the Hodge bundle (\cite{HaMo98} Formula~3.165) and a combination 
of the boundary
divisors. Since $C$ avoids both $H$ and the boundary,
this implies that for $f: X \to C$, the pullback of universal family
of curves to $C$, we have $\deg( f_* \omega_{X/C}) = 0$. Consequently,
the family of curves is isotrivial, contradicting the construction of $C$.
\par
Suppose that $Y_M \subset H$. Then $t^{[n]}$ is unramified
over $Y_M$. Since $Y_M$ differs from $Y$ only by the codimension
$\geq 2$-locus $W \cap Y$,  we have (choosing some basepoint $x$)
$$\pi_1(Y_M,x) \cong \pi_1(Y,t^{[n]}(x)) \cong \Gamma
\subset \SL(\fraco \oplus \fraco^\vee).$$
Hence $\pi_1(Y_M,x)$ is a lattice in a Lie group of rank $\geq 2$.
This contradicts a result of Farb and Masur (\cite{FaMa98}): The image
of such a lattice in the mapping class group has to be finite
i.e.\ after base change to a finite unramified cover of $M_3^{[n]}$
the Hilbert modular threefold $Y$ becomes simply connected. This
is absurd.
\par
The same argument works if the codimension of $H \cap Y_M$ in $Y_M$
is $\geq 2$. Consequently, there is a codimension one component of 
$H \cap Y_M$.
\end{proof}
\par
The second remark concerns the counterexample in \cite{McM03a}. There
McMullen found a curve $X_{\zeta_7}$ of genus $3$ with real multiplication,
in fact $X_{\zeta_7}: y^2=x^7-1$, and the following crucial relation between
the eigendifferentials $\omega_i=x^{i-1}dx/y$: We have
$$\omega_2^2 = \omega_1 \omega_3.$$
Such a relation was ruled out by Lemma~\ref{KSHM} over a
Hilbert modular variety. It thus can hold on $Y_M$ only
on the hyperelliptic locus. The condition in i) is a geometric
reformulation of this fact.
\par
Third, observe that even for $\fraco=\ZZ[\zeta_7+\zeta_7^{-1}]$ we cannot
disprove linearity using Theorem~\ref{obstacles}
for the part $B$ of the corresponding Hilbert modular threefold 
lying in the generic stratum. The problem is that 
$(X_{\zeta_7},\omega_2)$ lies in $\Omega M_3(2,2)$ rather than 
in the generic stratum. Of course, we can disprove linearity indirectly: If
$B$ was linear, it would be $\GL^+_2(\RR)$-invariant by 
Proposition~\ref{Kontconverse}, hence its
closure would be $\GL^+_2(\RR)$-invariant, too, contradicting
Theorem~7.5 in \cite{McM03a}. 
\par
\begin{proof}[Proof of Theorem~\ref{obstacles}]
Suppose that $B$ is linear.  Then by Theorem~\ref{DefofLLLemma} there 
is a linear subsystem $\LL_B$ of the VHS over $B$ such that $\psi$ maps
onto $\cOO_\PP(1)|_B \otimes \LL_B$. We write $t$ as shorthand for the
Torelli map $t^{[n]}$. By construction the 
fibers of $p:B \to Y_M \subset M_3$ are one-dimensional.
Since moreover,  via $\psi$, the relative tangent space maps into 
the $(1,0)$-part
of $\LL_B$, we conclude that $\LL_B$ has rank two and equals
the pullback to $B$ of, say, $t^* \LL_1$ on $Y$. Consequently
the fibers of $q$ have dimension $2$ when restricted to $B$.
Equivalently, the intersections of $B$ with the leaves of 
the foliation by relative periods is two-dimensional. This
excludes $B \subset \Omega M_3(4), \Omega M_3(2,2)$ 
and $\Omega M_3(3,1)$ and, using the full statement of
Theorem~\ref{DefofLLLemma}, we conclude that ii) holds.
\par
The second consequence of $\psi$ mapping surjectively to
$\cOO_\PP(1)|_B \otimes \LL_B$ is, using $\cOO_\PP(1)|_B \cong
p^*\cLL_1^{-1}|_B$, the surjectivity of 
$$ p^* T_{p(B)} \to R^1 f_* \omega_{X/B} \to 
p^*\cLL_1^{-1}|_B \otimes p^*\cLL_1^{-1}|_B\subset 
p^*\cLL_1^{-1}|_B \otimes R^1 f_* \cOO_X.$$ 
This map is a pullback of a map between bundles pulling back
from $Y_M$. Dualizing we obtain that the restriction to
$p(B)$ of the map $\psi^\vee_{Y_M}$ in the diagram
$$ 
\xymatrix{
t^*\cLL_1 \otimes t^*\cLL_1 \ar[r] \ar[drr]^>>>{\psi^\vee|_{Y_m}}
& t^*\cLL_1 \otimes f_* \omega_{X/Y_M}
\ar[r] \ar[d]& t^* \Omega^1_Y(\log S) \ar[d] \\
&  f_* \omega_{X/Y_M}^{\otimes 2} \ar[r] & \Omega^1_{Y_M}(\log t^{-1}(S))
\ar[d] \\
 & & \cOO_H \\
}
$$
is injective with locally free cokernel. The right column
is the short exact sequence coming from the ramification
along the hyperelliptic locus. The maps are the multiplication maps
and the Kodaira-Spencer maps by Lemma~\ref{givenbymult} and Lemma~\ref{KSHM}.
Since the composition map in the first row is injective and the 
image splits off by Lemma~\ref{KSHM}, the linearity condition implies i).
\par
The converse implication follows the same lines:
Condition i) ensures surjectivity of $\psi$ onto
$\cOO_\PP(1)|_B \otimes p^* \LL_1|_B$ and condition ii)
ensures that the image of $T_B \to T_S|_B$ is a
local subsystem mapping to $\cOO_\PP(1)|_B \otimes p^* \LL_1|_B$.
\end{proof}
\par
\begin{Rem}{\rm
Suppose that $\Omega M_3(\fraco) = \Omega M_3(1,1,1,1)$ for some order 
$\fraco$. Then the condition ii) in Theorem~\ref{obstacles}
may be rephrased in terms of variations of mixed Hodge structures.
Let $Y_{\rm gen} \subset Y_M$ be the open subset of the Hilbert
modular threefold in the image of the generic stratum, i.e.\ 
$$Y_{\rm gen} = p(B \cap \Omega M_3(1,1,1,1)),$$
where $B = B(\fraco)$ is the eigenform locus as above.
Let $f: X \to Y_{\rm gen}$ be the universal family and $D \subset X$
the universal divisor with the inclusion $j: (X \sms D) \to X$. 
The direct image  $R^1f_* (j_! \CC)|_{Y_{\rm gen}}$ carries a 
variation of mixed Hodge structures. While the weight one
quotient $R^1f_* \CC$ is known to split as in \eqref{HMVHS},
it is an interesting question to the determine direct summands of
$R^1f_* (j_! \CC)|_{Y_{\rm gen}}$.
\newline
In order to do so, one would like to know the fundamental group
of $Y_{\rm gen}$ as well as the representation
underlying $R^1f_* (j_! \CC)|_{Y_{\rm gen}}$.
\par
Although it does not seem directly relevant for linear manifolds, 
the same question arises for $g=2$: What is the fundamental group
of a Hilbert modular surface minus the locus of abelian surfaces
that decompose with polarization and minus the image of the locus
of eigenforms with a double zero. How can one describe its
representation underlying the mixed Hodge structure $R^1f_* (j_! \CC)|_{Y_M}$ ?
}\end{Rem}

\section{Linear manifolds not 
necessarily defined over $\RR$:
searching for a good definition}
\label{notoverRR}

In the previous sections we have tested manifolds to have
a linear structure defined over $\RR$ in the sense of 
Section~\ref{affstruc}. In the cases 
analyzed so far, such a manifold is all of a stratum, the
hyperelliptic locus or the linear structure is related to
uniformization. The purpose of this section is to make the
last sentence more precise in a way that linear structures
not necessarily defined over $\RR$ also fit into this context.
\par
{\bf Eigenforms over Ball quotients (\cite{DeMo86})}
\par
We consider the family of cyclic coverings of $\PP^1$ of degree $d$ 
ramified over the $N \geq 4$ points $\{\infty,0,1,x_1,\ldots, x_{N-3}\}$.
Hence let
$$M := \{(x_i)_{i=1}^{N-3} \,| \, x_i \neq x_j, x_i \not\in \{0,1,\infty\}\}  
\subset (\PP^1)^{(N-3)}$$
and let $f: X \to M$ be the universal family of such cyclic
coverings for a fixed type of covering. Here, the type contains
some more information besides the ramification orders. We
fix a type contained in the list on p.~86 of \cite{DeMo86} for $N \geq 5$
and recall the key properties of these families.
\par
In that situation $M$, or a completion of $M$ adding a finite number of 
boundary divisors, is shown by Deligne and
Mostow to be a quotient of the $(N-3)$-ball. In fact, the VHS over $M$
has a $\CC$-summand $\LL$ of type $(1,N-3)$ with the following
Higgs field $\Theta$. We let $\cLL = (\LL \otimes \cOO_M)^{(1,0)}_{\rm ext}$.
Then 
$$ \Theta: \cLL \cong \omega_X^{1/N}  \to 
(T_Y \otimes \omega_X^{1/N}) \otimes \Omega^1_X \cong (\LL)^{(0,1)} \otimes 
\Omega^1_X $$
the canonical diagonal embedding. The analog of the eigenform locus
in the Hilbert modular case is given by letting 
$\Omega M_g \supset B \to M$ be the total space of $\cLL$.
Thus $\dim B = N-2$.
\par
\begin{Prop} \label{DMislinear}
The total space $B$ is contained in some stratum $S \subset \Omega M_g$.
It is a linear submanifold of $S$ and {\em not} defined over $\RR$.
\end{Prop}
\par
\begin{proof} As in the case of Hilbert modular surfaces we have
$\cOO_\PP(1)|_B \cong p^* \cLL^{-1}$ and $\psi_\rel:T_{B/M} \to \cOO_\PP(1)|_B 
\otimes \cLL$ is the obvious isomorphism. We claim that $\psi$ descends
modulo $T_{B/M}$ to 
$$\ol{\psi}: p^* T_M \to \cOO_\PP(1)|_B \otimes p^*(\LL)^{(0,1)}.$$ 
Given the claim,
$\psi: T_B \to \cOO_\PP(1)|_B \otimes \LL$ is an isomorphism for
dimension reasons. The cokernel of 
$f^*\cLL \to \Omega^1_X$ is \'etale over $B$, in fact it coincides with
the preimage of the ramification points of the cyclic cover by construction.
Thus there is a unique local subsystem $\widetilde{\LL} 
\subset R^1 f_* j_!\, \CC$
that lifts $\LL$ to and through which $\psi$ factors. The linearity
of $B$ follows now from Theorem~\ref{DefofLLLemma}. 
\par
That $\LL$, hence $\widetilde{\LL}$, is not defined over $\RR$ follows 
immediately the fact that its $(1,0)$-part is of rank one, while
its complex conjugate is the $(0,1)$-part which is of rank $N-3$.
\par
It remains to check the claim. The quotient of $\psi$ on $p^* T_M$
is the Kodaira-Spencer map with the dual of the multiplication map
by Lemma~\ref{givenbymult} pulled back via $p^*$. It now suffices
to apply the description of 
the Higgs field $\Theta$ plus the relation between $\Theta$ and the 
Kodaira-Spencer map (see end of Section~\ref{affstruc}).
\end{proof}
\par
\begin{Example}
The bundle of one-forms generated by $\omega = ydx/(x(x-1)(x-\lambda)
(x-\mu)$ over the $2$-dimensional family of curves
$$ y^3=x^2(x-1)(x-\lambda)(x-\mu), \quad (\lambda,\mu) \in
  (\PP^1\sms\{0,1,\infty\})^2,\,\,\, \lambda \neq \mu $$
is a closed submanifold of $\Omega M_3(4)$ with linear structure. 
This structure is not defined over $\RR$, but only over $\QQ(\zeta_3)$.
\end{Example}
\par
\begin{Rem} \label{patho}
{\rm We now list phenomena and pathologies we do
not want a 'linear manifold' to have:
\begin{itemize}
\item[i)] For any flat surface, i.e.\ any pair $(X_0,\omega_0)$
of a smooth curve and a one-form, the set  $\CC^*\cdot(X_0,\omega_0)$
has a linear structure. But it is never defined over $\RR$ by Riemann's
period relations and its image in $M_g$ is a point.
\item[ii)] The construction of \cite{DeMo86} also works for several
types of $N=4$ and yields the following, see also \cite{BoMo05}.
For example, let $d=2mn$ for $m,n$ odd. Denote by $g: Y \to C$ the universal
family of cyclic coverings $C \to M$ pulled back from $M \cong 
\PP^1 \sms \{0,1,\infty\}$ to a finite unramified covering 
in order to kill non-unipotent monodromies.
From the construction as cyclic covering the VHS decomposes into
eigenspaces
$$ R^1 g_* \CC = \bigoplus_{(i,m)=1,(j,n)=1} \LL(i,j)$$
The $\LL(i,j)$ are rank two local systems, not defined over $\RR$
as local subsystems of $ R^1 g_* \CC$. In fact the $\LL(i,j)$
are defined over $\RR$ as abstract local system. A local
subsystem isomorphic to $\LL(i,j)$ and defined over $\RR$ sits
diagonally in a sum of $4$ copies of $\LL(i,j)$'s.
This is a key observation in \cite{BoMo05}, but not important
in this remark. 
\par
As a consequence, the $(1,0)$-part of $\LL(i,j)$ is one-dimensional
and we let $B(i,j)$ denote its total space.
From the construction as cyclic coverings it is immediate that
$B(i,j)$ lies completely in a stratum $S(i,j)$, which one depends on $(i,j)$.
For all pairs $(i,j)$ the total space is indeed a closed submanifold
of $S(i,j)$, since the fibers over $\{0,1,\infty\}$ are always singular
curves. Moreover, $S(i,j)$ carries a linear structure as we now check.
\par
The restriction $\psi_\rel$ of $\psi$ to $T_{B/p(B)}$ is an 
isomorphism, since $B$ is 
constructed as total space of a vector bundle on $p(B)=C$. Moreover it is 
well-known (see \cite{BoMo05} for a proof in this language) that
$$ \ol{\psi}: p^* T_{p(B)} \to \cOO_\PP(1)|_{B(i,j)} \otimes R^1 f_* \cOO_B$$
maps to, but not onto, $\cOO_\PP(1)|_{B(i,j)} \otimes \LL^{(0,1)}$.
For all $(i,j)$ the zeros of the differentials generating $\LL(i,j)^{(1,0)}$
are the preimages of $\{0,1,\infty,x_1\}$ via the cyclic covering.
In particular the difference of any two zeros is torsion in 
${\rm Pic}^0(Y/C)$. Hence $\LL(i,j)$ lifts to a local subsystem
of the local system of relative periods, as needed to apply
Theorem~\ref{DefofLLLemma}.
\par
But, fixing $d$, only for $4$ pairs $(i,j)$ the linear structure controls
the uniformization in the sense that the monodromy $\Gamma$ of $\LL(i,j)$ is
discrete and $\HH/\Gamma$ is a partial compactification
$C_0$ of $C$. The crux is that in all the other cases, 
the Kodaira-Spencer map $\ol{\psi}$ 
naturally extends to $C_0$, but it is no longer onto 
$\cOO_\PP(1)|_{B(i,j)} \otimes \LL^{(0,1)}$ at $C_0 \sms C$.
\item[iii)] It seems natural to remedy the above problem by
imposing the existence of a closure of $B$ with good properties.
Neither for $B$ (due to the $\CC^*$-action) nor for $p(B)$ (since 
all information about the one-forms is lost) we expect to have such
a compactification in general. On the other hand even in the intermediate 
case, just dividing by $\CC^*$, the quotient carries no longer a linear 
structure. But the tangent mappings still exist on $B/\CC^*$. This 
motivates the following definition, justified by Theorem~\ref{LINMF}. 
\end{itemize}
}
\end{Rem}
\par
Fix some stratum $S$ Let $\pi: S \to S/\CC^* \subset \PP \Omega M_g$ be
the quotient map by $\CC^*$ and let $\tau:  S/\CC^* \to M_g$ be
the forgetful map. We thus have a factorization $p = \tau \circ \pi$.
Recall that the local system $\LL_B$ depends only on the Hodge structure
of the fibers, not on the one-form chosen. It is thus a pullback
$\LL_B = p^* \LL$ of a local system $\LL$ on $p(B)$.
\par
\begin{Defi} \label{deflinmfk}
A submanifold $B$ in some stratum $S \subset \Omega M_g$ is called
a {\em linear manifold}, if $B$ is closed in $S$, if it inherits 
from $S$ a linear structure, if its image in $M_g$ is not reduced to a point 
and if there exists a compactification $Y$ of $\pi(B)$ with boundary 
$\Delta$ a
normal crossing divisor with the following property: Write $\Delta = T + U$
as disjoint union, such that the monodromy of $\LL$ is trivial 
around components of $T$ and unipotent around components of $U$.
There is a surjective map
$$\psi_Y: T_Y(-\log U) \to \cOO_{\PP \Omega M_g}(1)|_Y \otimes \tau^* (\LL),$$
whose restriction to $\pi(B)$ and pullback via $\pi$ yields the
quotient map of $\psi$
$$T_B/T_{B/\pi(B)} \to \cOO_{\PP}(1)|_B \otimes \LL_B/\psi(T_{B/\pi(B)}). $$
\end{Defi}
\par
\begin{Thm} \label{LINMF}
The following manifolds match the conditions of Definition~\ref{deflinmfk}
\begin{itemize} 
\item[i)] The connected components of strata and
the hyperelliptic locus are $\RR$-linear manifolds. 
\item[ii)] The eigenform locus over a Hilbert modular surface is a 
linear manifold defined over $\RR$.
\item[iii)] The canonical lift of a Teichm\"uller curve is a
linear manifold defined over $\RR$.
\end{itemize}
\end{Thm}
\par
\begin{proof}[Sketch of proof:]
In case i) for a stratum $S$, take the closure of $\pi(S)$ in the 
projectivized total space of the relative dualizing sheaf 
$\ol{\PP \Omega M_g}$ over the Deligne-Mumford compactification 
of $M_g$ and take a blowup of the lower-dimensional strata 
to make $\Delta$  a normal crossing divisor.
In that case, $\LL$ coincides with the full de~Rham cohomology.
Hence $U = \tau^{-1}(\partial M_g) \cap Y$. Around $U$ 
the Hodge metric has logarithmic growth. This implies that $\psi_Y$ is
surjective near $U$.
\par
For the hyperelliptic locus the same argument works using the
compactification by admissible double coverings (\cite{HaMo98}).
\par
For the Hilbert modular surfaces $B$ in ii) take $Y$ as in 
diagram~\eqref{levelstructure}
as compactification of $p(B)=Y_M$. Equation~\eqref{HilbVHS}
is precisely what we need.
\par
For iii) note that a Teichm\"uller curve $C \cong \pi(B) \cong p(B)$ 
does not intersect the
locus of curves with separating nodes. Consequently, 
$U = \ol{C} \sms C = \Delta$ and the desired surjection follows again from the
growth of the Hodge metric. 
\end{proof}
\par
To justify this complicated definition, recall that a Teichm\"uller
curve $C$ with compactification $\ol{C}$ and $\ol{C} \sms C = \Delta$
is characterized by having a rank two local system $\LL$ that is 
maximal Higgs, i.e.\ such that the Higgs field
$$ \theta: \LL^{(1,0)} \to \LL^{(0,1)} \otimes \Omega^1_{\ol{C}}(\log \Delta)$$
is an isomorphism. A $2$-dimensional linear manifold must be
{\em almost Teichm\"uller} in the following sense: i) By the second 
condition of Theorem~\ref{DefofLLLemma} and Definition~\ref{deflinmfk} it has
a rank two local system $\LL$ that is maximal Higgs with $U$, 
i.e.\ such that
$$ \LL^{(1,0)} \to \LL^{(0,1)} \otimes \Omega^1_{\ol{C}}(\log U)$$
is an isomorphism.
ii) It satisfies the first condition of Theorem~\ref{DefofLLLemma}
involving relative periods.
\par
We need to allow this modification of the usual definition, otherwise 
the locus of reducible Jacobians would cause a Hilbert modular surface
not to be a linear manifold.
From the construction as closed $\GL_2^+(\RR)$-orbits we deduce
that an $\RR$-linear almost Teichm\"uller curve is in fact
a Teichm\"uller curve. If we drop the condition linearity over
$\RR$ this is no longer the case, see Section~\ref{lowg2}.

\section{Linear manifolds in low 
dimension: $g=2$} \label{lowg2}

We show that the preparation made above yields a short
proof of the classification of manifolds with linear structure
in $\Omega M_2$ (originally due to McMullen, \cite{McM03b}) 
under the additional hypothesis that the manifold is algebraic
and with the slight generalization
of not restricting to linear manifolds defined over $\RR$.
\par
Referring to the notation of Theorem~\ref{DefofLLLemma}, we let 
for a linear manifold $B$ denote
$\tilde{d}:= {\rm rank}(\widetilde{\LL_B})$
and $d:={\rm rank}(\LL_B)$. If the linear structure is defined over
$\RR$, then $d$ must be even and $d/2$ is the dimension of the fibers
of $p: B \to p(B) \subset M_g$.
\par
\begin{Thm} \label{classg2}
A linear submanifold  $B$ in a stratum of $\Omega M_2$
is one of the following possibilities:
\begin{itemize}
\item[i)] One of the strata $\Omega M_2(2)$ or $\Omega M_2(1,1)$.
\item[ii)] The locus of eigenforms over a Hilbert modular surface
or over a surface parameterizing curves with a reducible Jacobian.
\item[iii)] The canonical lift of a Teichm\"uller curve to $\Omega M_2$
\item[iv)] A family of differentials over a Shimura curve parameterizing
Jacobians whose endomorphism ring is a quaternion algebra.
This quaternion algebra is indefinite or a matrix algebra.
\end{itemize} 
Precisely in the cases i), ii) and iii) the linear structure
is defined over $\RR$.
\end{Thm}
\par
An example of a quaternionic Shimura curve as in iv) is the
family $y^3=x(x-1)(x-\lambda)^2$. The total spaces of both eigenform
bundles for $\QQ(\zeta_3)$-multiplication 
generated by 
$$\omega_1 = ydx/x(x-1)(x-\lambda) \quad \text{and by} \quad
\omega_2 = y^2dx/x(x-1)(x-x_1)^2$$ 
are  $\QQ(\zeta_3)$-linear manifolds.
The zeros of $\omega_1$ are the preimages of $x=0$ and $x=1$, 
while the zeros of $\omega_2$ are the preimages of $x=\lambda$ and $x=\infty$. 
The quotient by the involution 
$$(x,y) \mapsto \left(\frac{\lambda}{x},\frac{(\lambda-x)(x-1)\lambda}{xy}
\right)$$
maps the family to a family of elliptic curves. Hence the  endomorphism 
ring is a matrix algebra. 
The proof of linearity is exactly the same as that of 
Proposition~\ref{DMislinear}.
\par
\begin{Q}
Is this the only example of such a curve in genus $2$? If
not, how to classify them? Are there examples where the
quaternion algebra is indefinite?
\end{Q}
\par
\begin{proof}[Proof of Theorem~\ref{classg2}] We list the cases 
for $(\tilde{d}, d)$: Case $d=1$ is excluded by the definition
of a linear manifold, compare Remark~\ref{patho}. If $d=3$, $B$
is not defined over $\RR$. The local system $\LL_B$ intersects its complex
conjugate non-trivially. Hence $B$ parameterizes a family of curves
whose Jacobian has a one-dimensional fixed part. For dimension
reasons, this is impossible both in  $\Omega M_2(1,1)$ and
$\Omega M_2(2)$.
\par 
In the stratum $\Omega M_2(2)$ the only remaining possibilities
are $(2,2)$ and $(4,4)$. Since $\Omega M_2(2)$ is irreducible,
$(4,4)$ corresponds to a stratum. The pair $(2,2)$ gives the
a $\CC^*$-bundle over a curve. We have to show that the
manifold is a Teichm\"uller curve or a quaternionic Shimura curve. 
This is the content of the Lemma below.
\par 
In the stratum $\Omega M_2(1,1)$ the possibilities are
first $(5,4)$, which is the whole stratum since it is irreducible,
second $(2,2)$, which is again a Teichm\"uller curve by the same arguments, 
and finally $(3,2)$. The case $(4,4)$ is impossible by the argument used in 
the proof of Theorem \ref{comphypend}: The fiber dimension
is $2$ in this case and the multiplication
map implies that $p(B)$ is at least $3$-dimensional.
\par
In the case $(3,2)$ the local system $\LL$
is irreducible for dimension reasons and since the
linear manifold $B$ does not map to a point in $M_2$. In particular,
the fibers of $B \to p(B)$ are generically one-dimensional.
Suppose that $\LL$ is not defined over $\RR$. Then $p(B)$
parameterizes curves whose Jacobian has endomorphisms
by a complex field. By \cite{Sh63} Proposition~19 this
endomorphism ring is in fact even larger. Such curves
are parameterized by Shimura curves. We conclude that
$p(B)$ is one-dimensional and obtain a contradiction
since $\dim(B)=3$.
\par 
Hence the family of Jacobians either splits or has 
real multiplication. In the second case, again for dimension reasons, 
$p(B)$ has to be dense in the
corresponding Hilbert modular surface. Recall that over a Hilbert
modular surface the VHS splits into $\LL_1 \oplus \LL_2$.
Hence for $i=1$  or $i=2$ the map
$$\psi_\rel: T_{B/p(B)} \to \cOO_\PP(1) \otimes p^*\LL_i^{(0,1)}$$
is an isomorphism. We deduce that $B$ is the total space
of the eigenform bundle over the Hilbert modular surface. 
The case of split Jacobians is similar.
\par
The linearity of the manifolds in the list has been shown in
Theorem~\ref{LINMF}.
\end{proof}
\par
\begin{Lemma} \label{TeichorQuat}
A linear manifold of dimension $2$, not necessarily $\RR$-linear, 
in a stratum of $\Omega M_2$ is the canonical  lift of a Teichm\"uller curve
or a quaternionic Shimura curve parameterizing
Jacobians whose endomorphism ring is a quaternion algebra, 
either indefinite or a matrix algebra.
\end{Lemma}
\par
\begin{proof}
Since multiplication by $\CC^*$ does not change the Hodge
structure, the local system $\LL_B$ (as in Theorem~\ref{DefofLLLemma})
is a pullback from $p(B)=:C$. Being isomorphic to $T_B$ and of
rank $2$, the  local system $\LL_B$ must be irreducible.
Hence either the VHS over $C$ splits over $\QQ$ and the Jacobians of
the universal family $X \to C$ split up to isogeny. 
Or the family of Jacobians of $f$ has multiplication by a field $K$
of degree $2$ over $\QQ$. 
\par
In the first case, though $\LL_B$ is defined over $\QQ$, we
cannot a priori not conclude that $\widetilde{\LL_B}$
is defined over $\RR$. But we know that the universal family
$f:X \to C$ admits a map $h:X \to E$ to a family of elliptic
curves $E \to C$. Moreover, since
$$\psi_\rel: T_{B/p(B)} \to \cOO_\PP(1) \otimes p^*\LL_B^{(0,1)}$$
is an isomorphism, we know moreover that the differentials 
parameterized by $B$ pull back from $E$. On fibers of $E$ there are 
only two absolute periods, $\RR$-linearly independent. Hence
if the relative period on $X$ depends linearly on the absolute
ones, it is possible to change the defining equation in order to obtain
linear dependence defined over $\RR$.
\par
The second case is $K$ is real, in particular $\LL_B$ is defined 
over $\RR$. Again we can not yet conclude that $\widetilde{\LL_B}$
is defined over $\RR$, too. But since
$$\psi_\rel: T_{B/p(B)} \to \cOO_\PP(1) \otimes p^*\LL_B^{(0,1)}$$
is an isomorphism, we know moreover that $B$ lies in the
eigenform locus of real multiplication. By Theorem~\ref{HM2linear} 
(or \cite{McM03b}) the eigenform locus is linear and defined over $\RR$. 
\par
The last case is $K$ not real. In this case the endomorphism ring of 
the family of Jacobians
is larger than $K$ (\cite{Sh63} Proposition~19). It is an indefinite
quaternion algebra, either a division algebra and the generic Jacobian
is simple or a matrix algebra and the family of Jacobians splits.
In both cases, abelian varieties with such an endomorphism ring
are parameterized by a countable union Shimura curves. Hence
for dimension reasons, $C$ coincides with such a Shimura curve. 
\end{proof}
\par

\section{Linear manifolds in low dimension: The hyperelliptic
locus $\mathcal{H}$ in the non-hyperelliptic component of $\Omega M_3(2,2)$}
\label{lininH22}

Combining the analysis of the VHS as in Section~\ref{lowg2} and
the use of the multiplication map as in Theorem~\ref{comphypend}, one 
can write 
down for a given stratum a list of possibilities of closed algebraic
submanifolds. We will do this for a special case, the hyperelliptic
locus $\mathcal{H}$ in the non-hyperelliptic component of $\Omega M_3(2,2)$, 
i.e.\ the one with odd spin structure $\Omega M_3(2,2)^{\rm odd}$. 
Recall that $\mathcal{H}$ is of dimension $5$ and
up to $\CC^*$ it is finite over the hyperelliptic locus. This choice
is motivated by \cite{HLM06}: On one hand, this locus is, besides $g=2$
probably the simplest one to classify $\GL_2^+(\RR)$-orbit closures
using connected sum constructions and Ratner's theorem for products of tori.
On the other hand, $\GL_2^+(\RR)$-invariant submanifolds of $\Omega M_2$
are classified, but $\GL_2^+(\RR)$-orbit closures in the
whole cotangent bundle to $M_2$, i.e.\ consisting of
a Riemann surface plus a quadratic differential, are not yet classified. 
By a well-known covering construction (recalled e.g.\ in \cite{HLM06}), 
they correspond bijectively to $\GL_2^+(\RR)$-orbit closures in 
 $\mathcal{H}$.
\par
It is convenient to introduce some more definitions. Given an
algebraic manifold $B$
with linear structure, consider the monodromy group $\Gamma$ of 
the local system $\LL_B$ as in Theorem~\ref{DefofLLLemma}. 
Recall from the proof of Lemma~\ref{notallthenEnd} that $\LL_B$ is isomorphic 
to a local system defined over a number field. In particular
$K:= \QQ(\tr(\gamma), \gamma \in \Gamma)$ is a number field, 
called the {\em trace field} of $B$. Here the more 
important notion is the field $F$ fixed by all elements $\sigma$ of 
${\rm Gal}(\ol{\QQ}/\QQ)$ such that $\sigma(\LL_B) \cong \LL_B$. 
We call $F$ the {\em fixed field} of $B$. Obviously $K \subset F$.
For Teichm\"uller curves (\cite{Mo06a} Lemma~2.3) and Hilbert modular
surfaces both fields coincide.
\par
\begin{Thm} \label{classH22}
A closed submanifold $B$ in $\mathcal{H}$ with linear
structure defined over $\RR$ is one of the following possibilities:
\begin{itemize}
\item[i)] the whole locus $\mathcal{H}$, or
\item[ii)] a Teichm\"uller curve, or
\item[iii)] a manifold of dimension $4$ with fixed field of
degree $1$ or $2$ over $\QQ$.
\end{itemize}
\end{Thm}
\par
\begin{proof} 
We first claim that the foliation by absolute periods and
the hyperelliptic locus intersect transversely in 
$\Omega M_3(2,2)^{\rm odd}$. Consider the diagram~\eqref{dualidentification}
restricted to a point $(X_0: y^2=\prod (x-x_i), \omega_0)$ in 
$\Omega M_3(2,2)^{\rm odd}$. The hyperelliptic
involution $h$ acts on all bundles involved. Since $D-D_\red$ is
fixed by $h$, the involution $h$ fixes $\Ker(\pi)$, too. 
On the whole de~Rham cohomology
$h$ acts by $(-1)$, since $x^{i-1} dx/y$ for $i=1,2,3$ is a basis of 
the one-forms and $[y/x^i] \in \Gamma(X_0 \sms \{0,\infty\},\cOO_{X_0})$ 
for $i=1,2,3$ is a  basis, using Czech cohomology, 
of $H^1(X_0,\cOO_{X_0})$. Hence the
image of $q^\vee|_{\Ker(\pi)}$ is in the $+1$-eigenspace of
$H^0(X_0, (\Omega^1_{X_0})^{\otimes 2})$. This eigenspace is 
precisely the cotangent space to the hyperelliptic locus (e.g.\ \cite{OS80}). 
Since, dually, the kernel of $q$ defines the tangent space to the
foliation by absolute periods, the claim follows.
\par
From the claim we conclude that the possible dimensions of
$\widetilde{\LL}_B$ and $\LL_B$ are $\tilde{d}=d \in \{6,2,4\}$.
These give the rough classification.
\par
It remains to check the claim on the degree of the fixed field $F$
in the case $d=4$. If $[F:\QQ]=3$, then curves parameterized by $B$
 have real multiplication by $F$. For dimension reasons $p(B)$
equals a dense subset of a Hilbert modular threefold. This contradicts
Proposition~\ref{hypintHM3}.
\end{proof}
\par
An example for case iii) with fixed field $\QQ$ is to take for $p(B)$ 
a component of the Hurwitz space of unramified double coverings of a 
genus~$2$ surface and for $B$ the pullback of elements in
$\Omega M_2(1,1)$ under this covering map.
For a suitable choice of the double covering, the surfaces in
$\Omega M_3(2,2)$ thus obtained are indeed in $\Omega M_2(2,2)^{\rm odd}$
and hyperelliptic.
\par
\begin{Cor} \label{pAclosure}
Let $\Delta$ be a Teichm\"uller disc 
generated by a pair $(X_0,\omega_0)  \in \mathcal{H}$ which
is stabilized by a pseudo-Anosov diffeomorphism with trace field
of degree $3$ over $\QQ$. If the closure of the orbit 
$\GL_2^+(\RR) \cdot (X_0,\omega_0)$ is an algebraic manifold,
then it is either the canonical lift of a Teichm\"uller curve to $\mathcal{H}$
or as big as possible, i.e.\ $\ol{\GL_2^+(\RR) \cdot (X_0,\omega_0)} = 
\mathcal{H}$.
\end{Cor}
\par
\begin{proof} Suppose the claim is wrong. Then by Theorem~\ref{classH22}
we deduce that the closure $B:= \ol{\GL_2^+(\RR) \cdot (X_0,\omega_0)}$
has dimension $4$ and $\dim p(B) =3$. 
The pseudo-Anosov diffeomorphism corresponds to a closed geodesic
$\gamma$ in the image of $\Delta$ of $\Delta$. The monodromy of
the subspace $V:=\langle {\rm Re} \omega_0, {\rm Im} \omega_0 
\rangle$ along $\gamma$ has a trace $\tr(\gamma)$ that
generates a field with $[\QQ(\tr(\gamma)):\QQ]=3$. Suppose that for
$\sigma \in {\rm Gal}(\ol{\QQ}/\QQ)$ we had $\sigma(\LL_B) \cong \LL_B$.
Then this holds in particular when restricting the local system
to $\gamma$. Moreover it holds for the invariant subspace $V$
of $\LL_B|_{\gamma}$. Hence $\sigma$ fixes $\QQ(\tr(\gamma))$ and
statement iii) in Theorem~\ref{classH22} yields the contradiction.
\end{proof}
\par
In \cite{HLM06} this closure statement was shown for a 
particular Teichm\"uller disc, that does not descend to a
Teichm\"uller curve, without any a priori manifold hypothesis 
on the closure.
\par

\par
Martin M{\"o}ller: \newline
Universit{\"a}t Duisburg-Essen, \newline 
FB 6 (Mathematik) \newline 
Campus Essen \newline
45117 Essen, Germany \newline
\par
Current address: \newline
Max-Planck-Institut f\"ur Mathematik \newline
Postfach 7280 \newline
53072 Bonn, Germany \newline
email: moeller@mpim-bonn.mpg.de
\end{document}